%% file: hpc_darcy_mrcm_V1.tex
\documentclass[a4paper,english]{article}
\usepackage[a4paper,left=1in,right=1in,top=1in,bottom=1in]{geometry}
\usepackage{graphicx}
\usepackage{amsmath,amsthm,amsfonts,amssymb}

\usepackage{multirow}
\usepackage{url}
\usepackage{booktabs}
\usepackage{algorithm2e}

\usepackage{color}

\newenvironment{nalign*}{
	\begin{equation*}
	\begin{aligned}
}{
	\end{aligned}
	\end{equation*}
	\ignorespacesafterend
}

\newcommand{\parder}[2]{\frac{\partial #1}{\partial #2}}

\newcommand{\Sw}{S_\text{w}}

\newtheorem{remark}{Remark}

\numberwithin{equation}{section}

\graphicspath{{figs/}}

\title{
  {Towards HPC simulations of Billion-cell Reservoirs by Multiscale Mixed Methods} \\
}

\author{Alfredo Jaramillo$^{a}\footnote{Alfredo Jaramillo. ajaramillopalma@gmail.com. Institute of Mathematics and Computer Science of the University of São Paulo. Av. Trab. São Carlense, 400 - Centro, São Carlos - SP, Brazil, 13566-590.}$, Rafael T. Guiraldello$^a$, Stevens Paz$^a$, Roberto F. Ausas$^a$, \\Fabricio S. Sousa$^a$, Felipe Pereira$^b$, Gustavo C. Buscaglia$^a$\\
{\small
(a) Instituto de Ci\^encias Matem\'aticas e de Computa\c c\~ao, Universidade de S\~ao Paulo
}\\
{\small
(b) Department of Mathematical Sciences, The University of Texas at Dallas 
}
}

\begin{document}
\vspace{3cm}

\maketitle


\begin{abstract}

  A three dimensional parallel implementation of Multiscale Mixed Methods
  based on non-overlapping domain decomposition techniques is proposed for multi-core computers and
  its computational performance is assessed by means
  of numerical experimentation. As a prototypical method, from which many
  others can be derived, the Multiscale Robin Coupled Method is chosen and its implementation
  explained in detail. 
  Numerical results for problems ranging from millions up to more than 2 billion computational cells
  in highly heterogeneous anisotropic rock formations based on the SPE10 benchmark are shown.
  The proposed implementation relies on direct solvers for both
  local problems and the interface coupling system. We find good 
  weak and strong scalalability as compared against a state-of-the-art global
  fine grid solver based on Algebric Multigrid preconditioning in
  single and two-phase flow problems. 
  
\medskip


\end{abstract}

\section{Introduction}\label{sec:introduction}

Reservoir simulation is a challenging problem that demands intensive use of computational resources.
For very thick reservoirs, as is the case in the Brazilian pre-salt layer, discrete models
require a number of computational cells in the order of billions so as to achieve results
with the required level of accuracy. Numerical methods based on a single global scale discretization
along with preconditioned Krylov iterative strategies to solve the resulting system of equations may
fail to handle such extreme levels of dicretization efficiently.

In the last decade, multilevel and multiscale methodologies based on domain
decomposition have provided an alternative with potential to deliver
reasonably accurate results \cite{manea2016,durlofsky2007,lie2017}
with good scalability properties. Multiscale methods decompose the global fine-scale
problem into several local problems, whose solutions
are called {\it Multiscale Basis Functions} (MBFs) and are computed by setting suitable boundary conditions
on the local domains. The MBFs computations are independent from each other and therefore its parallelization is
straightforward in multi-core computers. On the interface between the local domains the compatibility conditions, namely,
pressure and/or flux continuity, are imposed weakly on a coarser scale related to that of the domain
decomposition. The resulting relaxed interface problem is solved using the computed MBFs and its
coarse-grid interface solution is used to locally reconstruct the multiscale solution. 
Multiscale methods based on domain decomposition techniques aiming at the accurate approximation of velocity and pressure fields have been
the
focus of intensive research over the last years \cite{Pereira::2014,guiraldello2018multiscale,arbogast2007,JENNY2003,EFENDIEV2013,HOU1997,harder2013}.
However, in the context of reservoir simulation
involving highly anysotropic rock formations exhibiting regions with multiple channels and obstacles, 
the accuracy of the multiscale methods can seriously deteriorate \cite{kippe2008, rocha2020}.
Also, a lack of High Performance Computing reports involving high-resolution problems
(in the order of one billion cells) can be identified in the literature where most
of the reported results are only two dimensional and/or involve small/middle size problems \cite{manea::2016,manea::2017,Puscas::2018}.
The behavior in terms of accuracy and scalability of such methods in large scale complex three
dimensional formations, which are the target problems these methods have originally been devised to,
is still an area of active research \cite{abreu2020recursive}.

We focus our analysis on the Multiscale Robin Coupled Method (MRCM)
\cite{guiraldello2018multiscale}, a recently proposed method that generalizes the
Multiscale Mixed Method (MuMM) \cite{Pereira::2014}. The compatibility conditions on the interface
are enforced by means of a Robin type boundary condition which depends on an algorithmic parameter that can
be adjusted locally or globally, leading to a family of methods of which the Multiscale Mortar Mixed
Finite Element Method \cite{arbogast2007} and the Multiscale Hybrid-Mixed Finite Element Method \cite{harder2013}
are particular cases. This feature of the method is attractive, since it provides a general framework
that can be useful to study the accuracy and performance of multiscale mixed methods based on
non-overlapping domain decomposition.

Several aspects of this method have been explored in the last few years. For instance, the choice of the MRCM algorithmic
parameter has in part been surveyed in \cite{guiraldello2018multiscale}. Also, Rocha et al \cite{rocha2020} have
proposed an adaptive strategy to improve accuracy by setting this algorithmic parameter as a function of
the permeability contrast.
Concerning the interface spaces for the flux and pressure coupling unknowns on the skeleton
of the domain decomposition, Guiraldello et al \cite{guiraldello2019} and
Rocha et al \cite{rochaphd,rocha2021prerint} have introduced, respectively, informed and physically based
basis functions that can yield more accurate results than polynomial spaces.
Finally, the coupling to a transport solver in the linear passive case \cite{guiraldello2020}
and in the non-linear (two phase flow) case \cite{rocha2020} have also been carefully investigated.

The aforementioned contributions illustrate the potential of multiscale mixed approximations for the flow calculations
in high-contrast permeability fields. However, in this papers the authors assess the method by looking at accuracy 
and convergence to fine grid solutions in two dimensional setting involving no more than $\sim10^4$ unknowns,
the three dimensional case being left aside. In this work we fill in this gap by showing high-resolution
three dimensional numerical results produced by a HPC implementation of the MRCM.
Weak and strong scalability as well as accuracy of the method are investigated by comparison
to fine grid solutions. These solutions are based on a classical finite volume discretization scheme
considering sizes ranging from $\sim 20$ million to a few billion
computational cells, defined by suitably refining the challenging SPE10 industry benchmark \cite{christie2001}.

The structure of this paper is as follows: Section \ref{sec:mathematical_models} provides the mathematical models
and numerical techniques used in this work, which involves a global elliptic solver (the Fine Grid method) and
a three-dimensional implementation of the MRCM. Details on several implementation aspects are provided.
Computational times and efficiency are reported in Section \ref{sec:hpc-results-darcy} for both elliptic solvers.
Special attention is given to the weak and strong scalability and the gain factor obtained in solving problems by a multiscale method.
To assess the precision of the multiscale solution two-phase flow simulations are also performed within an IMPES \cite{Sheldon1959,IMPES,chen2004improved,PAZ2021} approximation 
and the results for each one of the solvers are compared.
Finally, some concluding remarks are drawn in Section \ref{sec:conclusions}.

\section{Mathematical model and numerical schemes}
\label{sec:mathematical_models}
\input{tex/hpc_math_num_models}

\section{HPC results}\label{sec:hpc-results-darcy}

The results reported in this work were obtained in the \emph{Euler Supercomputer}\footnote{\url{https://sites.google.com/site/clustercemeai}},
that is located at the Institute of Mathematical and Computer Science (ICMC) of the University of
S\~ao Paulo at S\~ao Carlos\footnote{\url{https://www.icmc.usp.br}} The hardware specifications of listed in Table \ref{tab:computational_nodes}.
\begin{table}[h!]
	\begin{center}
		\begin{tabular}{ll}
			\toprule
			\multicolumn{2}{c}{\emph{Hardware specifications} }          \\ 
			\midrule
			Socket model & Intel Xeon \verb|E5-2680v2|    \\
			Clock frequence & $2.8$ GHz  \\
			Cores per socket & 10     \\
			Sockets per node & 2     \\
			RAM per node & 128 Gbytes  \\
                        Communications & InfiniBand \\
                        \bottomrule
		\end{tabular}
		\caption{Specifications of the computational nodes were the implementations were tested.}\label{tab:computational_nodes}
	\end{center}
\end{table}

\subsection{Numerical setup}

\begin{figure}[h]
	\includegraphics[width=\textwidth]{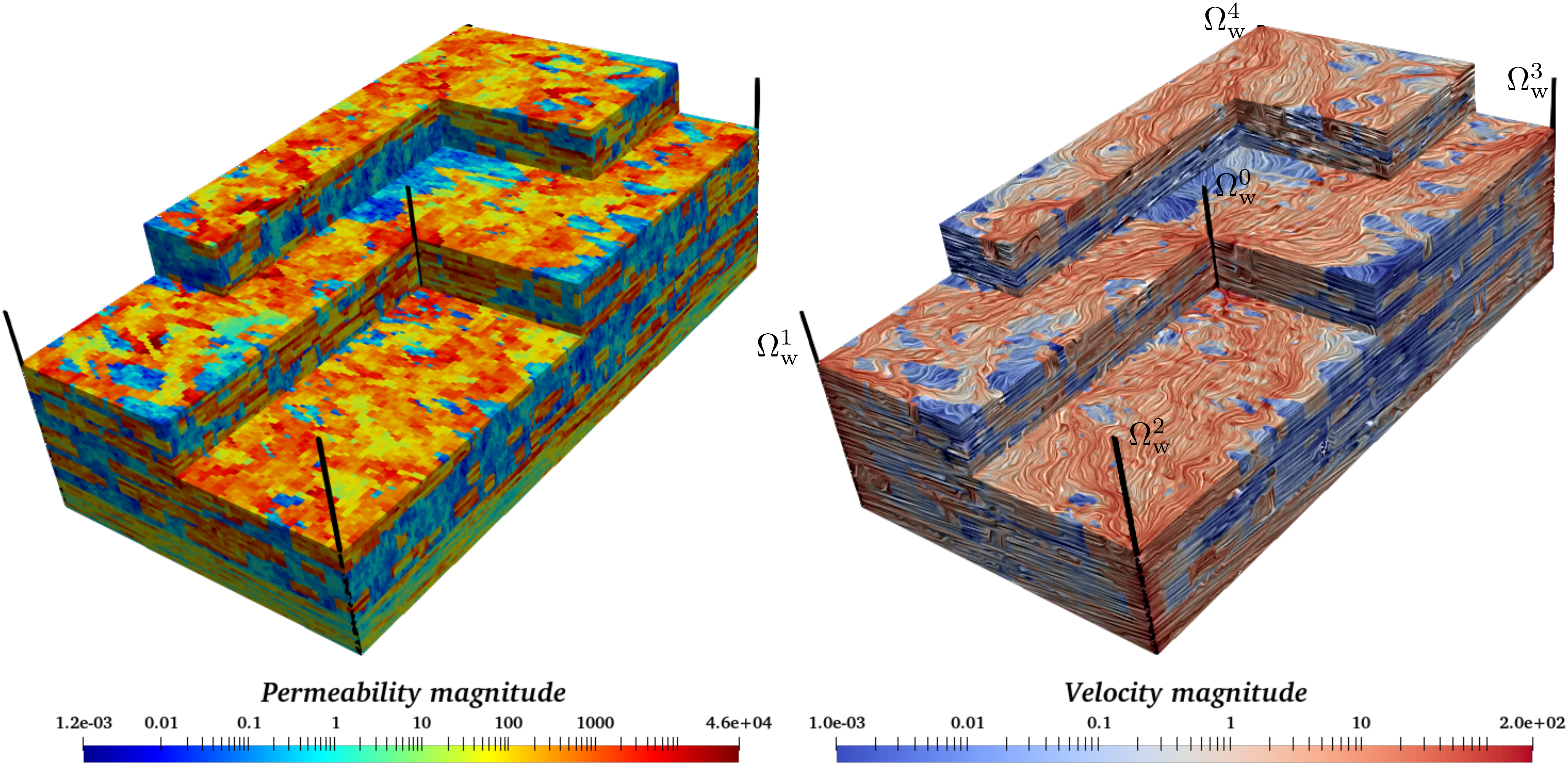}
	\caption{Permeability (magnitude) field $\mathbf{K}$ cuts and wells disposition (left). Resulting flux magnitude for the MRCM (right).}\label{fig:spe10_wells}
\end{figure}

The numerical experiments are performed on a 5-spot geometry posed on a parallelepiped region $\Omega = [0,1200]\times[0,2200]\times[0,120]$ (ft)
with absolute permeability $\mathbf{K}$ taken as the diagonal tensor provided by the SPE10 project \cite{christie2001} posed on a reference grid
of $60\times220\times60$. Instead of considering the whole set of layers, for convenience in the design of the numerical
experiments, only layers 26 to 85 have been selected in the $x_3$ direction.
This field provides very heterogeneous channelized structures with high permeability contrast and is a standard benchmark for numerical simulations in reservoir engineering.
In order to conduct numerical experiments on finer grids than the original $60\times220\times60$, an $L^2$ projection of the reference permeability is performed to populate
the tensor $\mathbf{K}$.
Five vertical well columns modeled as volumetric sources/sinks, corresponding to the rigth hand side $f$ in \eqref{eq:sec2:continuity}, with
dimensions $20\times10\times120$ (ft) are distributed in the physical domain as follows: one central well $\Omega_{\text{w}}^{0}$ injects water
(this is $S_{w}(\mathbf{x},t) = 1$, for $\mathbf{x}\in\Omega_{\text{w}}^{0}$, $t\geq 0$) such that one pore volume is injected (PVI)
every five years and four oil producer wells $\Omega_{\text{w}}^{i}$, $i=1,2,3,4$ located at the corners of the domain as shown
in Fig \ref{fig:spe10_wells}. No-flow boundary conditions are imposed on the external walls and the initial condition is set to $S_{w}(\mathbf{x},0) = 0$
on $\Omega \setminus \Omega_{\text{w}}^{0}$.

As for the MRCM, the interface spaces $\mathcal{P}_H$ and $\mathcal{U}_H$ are spanned by constant functions defined over
the domain decomposition with support on $\bar{H}_w\times\bar{H}_{\bar{w}}$, $w,\bar{w}=x_1,x_2,x_3$.
One can check \cite{guiraldello2018multiscale,guiraldello2019,rochaphd} for alternative choices of interface spaces.
In line with previous works, the Robin condition parameter is defined locally over the interface grid as
\begin{equation}
\mathbf{\beta}_l(\mathbf{x}) = \dfrac{\alpha\,H_l}{K_l(\mathbf{x})},
\end{equation}\label{Robin_par}
where $\alpha$ is an algorithmic parameter, that unless otherwise mentioned, is set to
$1.0$, $H_l =\hat{\mathbf{n}}^\intercal\cdot\left[H_{x_1} \ H_{x_2} \ H_{x_3}\right]^\intercal$
and $K_l(\mathbf{x})=\hat{\mathbf{n}}^\intercal\cdot\mathbf{K}(\mathbf{x})\cdot\hat{\mathbf{n}}$.

Accuracy of the multiscale solutions is assessed by comparing them against the fine grid solution in terms of relative errors
measured in the standard $L^2$-norm for pressure and the weighted $L^2$-norm given by $\lVert{\bf v}\rVert_{K^{-1},\omega} = \left( \int_\omega K^{-1}|{\bf v}|^2\right)^\frac{1}{2}$
for velocity.

\subsection{Darcy solver: Weak and strong scaling}

\subsubsection{Weak scaling}
\label{sec:weak_scaling}

The weak scaling properties of the algorithms is assessed in three different problem configurations
as indicated in Table \ref{wsc_tab}. The number of cores is increased as the fine grid is refined
in the same proportion, so perfect weak scalability would imply the computational time to remain constant.
The processor decomposition indicates the number of
MPI processes and the core decomposition refers to the number of subdomains solved per core.
{\color{black} The global number of unknowns for each case and the size of the local problems
are also indicated in the table}. 
\begin{table}[h!!!]
	\begin{center}
		\begin{tabular}{ccrrr}
			\toprule
			& core decomp. & 4$\times$2$\times$4 &8$\times$2$\times$4&8$\times$2$\times$8\\ \midrule
			\# cores & processor decomp. & $P^{wsc}_1$  & $P^{wsc}_2$ & $P^{wsc}_3$                            \\ \midrule
			32    & 2$\times$4$\times$4    & {\small  80$\times$1760$\times$160({\bf 22.5M})}  & {\small 160$\times$1760$\times$160({\bf 45M})}   & {\small 160$\times$1760$\times$320({\bf 90M})}    \\
			64    & 4$\times$4$\times$4    & {\small 160$\times$1760$\times$160({\bf45M})}  & {\small 320$\times$1760$\times$160({\bf 90M})}   & {\small 320$\times$1760$\times$320({\bf 180M})}    \\
			144   & 6$\times$4$\times$6    & {\small 240$\times$1760$\times$240({\bf 101M})} & {\small 480$\times$1760$\times$240({\bf 203M})}  & {\small 480$\times$1760$\times$480({\bf 406M})}   \\
			288   & 6$\times$4$\times$12   & {\small 240$\times$1760$\times$480({\bf 203M})} & {\small 480$\times$1760$\times$480({\bf 406M})}  &{ \small 480$\times$1760$\times$960({\bf 811M})}  \\
			576   & 8$\times$6$\times$12  & {\small 320$\times$2640$\times$480({\bf 406M})} & {\small 640$\times$2640$\times$480({\bf 811M})} & {\small 640$\times$2640$\times$960({\bf 1620M})}\\
			960  & 8$\times$10$\times$12 & {\small 320$\times$4400$\times$480({\bf 676M})} & {\small 640$\times$4400$\times$480({\bf 1350M})} & {\small 640$\times$4400$\times$960({\bf 2700M})} \\ \bottomrule
		\end{tabular}
		\caption{Domain decomposition and global grid for three problems for a weak scaling analysis.}
		\label{wsc_tab}
	\end{center}
\end{table}

The average computational times of the different stages are reported in Table \ref{wsc_res}. In this table we show the total times for both fine grid and
multiscale solver (MRCM), the time spent in MBFs computation and for the interface linear system solution.
These times are obtained over several simulations such that deviations
from the mean are negligible. The relative pressure and velocity errors $\|e_p\|_{L^2(\Omega)}$ and
$\|e_{\bf v}\|_{K^{-1},\Omega}$ are also reported.

\begin{table}[h]
	\centering
	\begin{tabular}{ccccccccrr}
		\toprule
		& {\bf \#cores} & {\bf Fine } & {\bf MRCM} & {\bf MBFs} & {\bf Int.} & $\|e_p\|_{L^2(\Omega)}$ & $\|e_{\bf v}\|_{K^{-1},\Omega}$ & $N_{\Omega}$ & $N_\Gamma$ \\ 
		\midrule
		\multicolumn{1}{c}{$P^{wsc}_1$} &         &   {\small $\frac{\sigma}{\mu}=0.33$}    &   {\small $\frac{\sigma}{\mu}=0.14$}   & {\small $\frac{\sigma}{\mu}=0.03$}     &      &          &          &        &    \\ \cline{1-1}
		& 32      & 14.4  & 7.2  & 6.7  & 0.0  & 0.45 & 0.38 & {\bf 22.5M}  & {\bf 5.5K}  \\
		& 64      & 17.3  & 7.1  & 6.5  & 0.1  & 0.63 & 0.50 & {\bf 45.1M}  & {\bf 11.3K}  \\
		& 144     & 22.0  & 7.4  & 6.6  & 0.3  & 0.60 & 0.54 & {\bf 101M}   & {\bf 25.8K} \\
		& 288     & 27.4  & 7.9  & 6.8  & 0.6  & 0.54 & 0.54 & {\bf 203M}  & {\bf 51.8K}  \\
		& 576     & 28.6  & 8.9  & 7.0  & 1.2  & 0.51 & 0.53 & {\bf 406M}  & {\bf 105.6K}  \\
		& 960     & 38.9  & 10.2 & 7.1  & 2.5  & 0.37 & 0.50 & {\bf 676M}   & {\bf 178.0K} \\ \cline{2-10} 
		\multicolumn{1}{c}{$P^{wsc}_2$} &         &   {\small $\frac{\sigma}{\mu}=0.25$ }   &   {\small $\frac{\sigma}{\mu}=0.18$}   &  {\small $\frac{\sigma}{\mu}=0.03$}    &      &          &          &        &    \\ \cline{1-1}
		& 32      & 35.9  & 14.2 & 13.0 & 0.1  & 0.20 & 0.44 & {\bf 45.1M} &  {\bf 11.3K} \\
		& 64      & 50.4  & 14.1 & 12.9 & 0.2  & 0.74 & 0.57 & {\bf 90.1M} & {\bf 22.8K} \\
		& 144     & 67.0  & 14.6 & 13.0 & 0.5  & 0.77 & 0.63 & {\bf 203M} & {\bf 51.8K}  \\
		& 288     & 83.3  & 15.6 & 13.4 & 1.0  & 0.73 & 0.63 & {\bf 406M}  & {\bf 104.5K}  \\
		& 576     & 73.3  & 17.7 & 13.8 & 2.6  & 0.70 & 0.65 & {\bf 811M} & {\bf 212.3K}  \\
		& 960     & 68.9  & 22.6 & 14.1 & 7.2  & 0.57 & 0.62 & {\bf 1350M} & {\bf 358.0K} \\ \cline{2-10} 
		\multicolumn{1}{c}{$P^{wsc}_3$} &         &  {\small $\frac{\sigma}{\mu}=0.27$}     &   {\small $\frac{\sigma}{\mu}=0.21$}   & {\small $\frac{\sigma}{\mu}=0.03$}      &      &          &          &          & \\ \cline{1-1}
		& 32      & 82.1  & 28.1 & 25.6 & 0.3  & 0.50 & 0.43 & {\bf 90.1M}  & {\bf 22.8K} \\
		& 64      & 112.2 & 27.8 & 25.4 & 0.5  & 0.67 & 0.56 & {\bf 180M}   & {\bf 46.1K} \\
		& 144     & 168.5 & 28.6 & 25.5 & 0.8  & 0.73 & 0.63 & {\bf 406M}  & {\bf 104.5K}  \\
		& 288     & 196.7 & 30.5 & 26.4 & 1.7  & 0.72 & 0.63 & {\bf 811M}   & {\bf 209.7K} \\
		& 576     & 179.7 & 35.8 & 27.1 & 5.9  & 0.69 & 0.65 & {\bf 1620M}  & {\bf 426.6K} \\
		& 960     & 168.7 & 47.1 & 27.6 & 16.8 & 1.07 & 0.61 & {\bf 2700M}  & {\bf 718.6K} \\ 
		\bottomrule 
	\end{tabular}
	\caption{Computational times for the fine grid solver, the MRCM, the time to build
		the MBFs and the time for the interface solver as well the relative error for pressure
		and velocity.}\label{wsc_res}
\end{table}

As one can notice the multiscale method has a better behavior with respect to weak scalability when
compared to the fine grid solver. Computational times for both methods increase
with the number of cores. The weak scalability with more than $288$ cores for the MRCM declines.
From Table \ref{wsc_res} we also conclude that the total time taken by the MRCM is mainly dominated
by the MBFs computation and the interface solver. As a measure of scalability performance,
we report for each case the quantity $\sigma/\mu$, which represents the relative deviation
from the mean time along the scalability test.
On one hand, the MBFs computation exhibits excellent weak scalability,
an expected behavior as the number of local unknowns remains constant
in each core. On the other hand, the interface resolution exhibits an inferior weak scalability behavior as the size
of the linear system increases with the number of subdomains to be coupled in the multiscale
process. This is explained by looking at the interface computing times. We observe a linear behavior with respect
to $N_{\Gamma}$ for $N_{\Gamma} \lesssim 200$K and a transition to quadratic behavior from that point
as previously shown in section \ref{sec:msc_solver_detail}. This has a deleterious effect for runs above $\sim 800$M.
One should also notice that for the three problems in Table \ref{wsc_res} the subdomains have the same $22$K unknowns,
but the cores handle twice the number of subdomains for $P^{wsc}_j$ when compared to $P^{wsc}_{j-1}$.

Next, we take problem $P^{wsc}_2$ and test two additional Darcy decompositions per processor, namely, $8\times1\times4$
and $8\times4\times4$ that lead to local Darcy problems of $44$K and $11$K unknowns, besides the original
setting of $8\times2\times4$ of $22$K.
Table \ref{wsc_p2} shows the total time for the fine grid, MRCM and time spent in the MBFs and interface solver phases.
Note that the MBFs computations is around $14$ seconds, since for these
sizes of local problems, we have not entered into the quadratic behavior observed in Figure \ref{fig:sec2:times-mumps}.
Taking larger subdomains ($\gtrsim 44$K) implies in an increased MBFs computing time. This was actually
confirmed by additional experiments, not shown here for the sake of brevity.
Moreover, the time for the interface solver gets worse with the increase in the number of subdomains, as expected.

\begin{table}[h]
	\centering
	\begin{tabular}{cc|ccc|ccc|ccc}
		\toprule
		\multicolumn{2}{c}{$P^{wsc}_2$} & \multicolumn{3}{c}{$8\times1\times4~($44$\mbox{K})$} & \multicolumn{3}{c}{$8\times2\times4~($22$\mbox{K})$} & \multicolumn{3}{c}{$8\times4\times4~($11$\mbox{K})$} \\ \midrule
		{\bf \#cores}         & {\bf Fine}         & {\bf MRCM}      & {\bf MBFs}      & {\bf Int.}     & {\bf MRCM}     & {\bf MBFs}      & {\bf Int.}     & {\bf MRCM}      & {\bf MBFs}   & {\bf Int.}     \\ \midrule
		32              & 35.9          & 15.2      & 14.0      & 0.0           & 14.2      & 13.0      & 0.1           & 13.9      & 12.6      & 0.3           \\
		64              & 50.4          & 15.2      & 13.8      & 0.1           & 14.1      & 12.9      & 0.2           & 14.0      & 12.5      & 0.5           \\
		144             & 67.0          & 15.4      & 14.0      & 0.2           & 14.6      & 13.0      & 0.5           & 15.1      & 12.7      & 1.3           \\
		288             & 83.3          & 16.4      & 14.7      & 0.4           & 15.6      & 13.4      & 1.0           & 16.9      & 13.0      & 2.8           \\
		576             & 73.3          & 17.5      & 15.1      & 0.9           & 17.7      & 13.8      & 2.6           & 25.3      & 13.3      & 10.7          \\
		960             & 68.9          & 18.9      & 15.5      & 1.8           & 22.6      & 14.1      & 7.2           & 40.8      & 13.6      & 26.0          \\ \bottomrule
	\end{tabular}
	\caption{Computational times to build the MBFs and solve the interface system for different local core decomposition of problem $P^{wsc}_2$.
          In parenthesis we show the number of unknowns of each subdomain. }\label{wsc_p2}
\end{table}

\subsubsection{Strong scaling}

{\color{black}
For the strong scaling assessment we consider three different global problems and two coarse scale decompositions.
For a given coarse decomposition, the number of interface unknowns $N_{\Gamma}$ and the size of the
local problems $N_{\Omega}$ are fixed. These are displayed in Table \ref{strong_tab}.
Strong scaling results for the different cases are summarized in Table \ref{sts}.
First, the optimum linear time decay is never observed for the three problems under consideration,
neither for the fine grid nor for the multiscale solver. Nonetheless, the situation is significantly
better with the multiscale solver as the total computing time goes as $\sim n^{-0.5}$ for problem
$P_1^{ssc}$ and $\sim n^{-0.85}$ for problems $P_2^{ssc}$ and $P_3^{ssc}$, $n$ being
the number of cores. For the fine grid solver the computing time goes as $n^{-0.3}$ for
$P_1^{ssc}$ and $n^{-0.6}$ for $P_2^{ssc}$ and $P_3^{ssc}$.

Note that, for a given number of cores, if we compare $P^{ssc}_j$ with $P^{ssc}_{j-1}$
an advantage of the multiscale solver with respect to the fine grid is observed
in terms of scalability. We remark that the first point of $P_3^{ssc}$ was omitted
because it did not even fit in memory.

When fixing the problem size $P_j^{ssc}$, there is also a clear improvement in computational time
of the multiscale solver as indicated in the Msc-Gain column, that is the ratio
of the fine grid to the multiscale total times.
Notice that the time spent in the interface solver is fixed for each coarse decomposition. This
time is given in Table \ref{sts}. When the processor grid increases each core ends up handling
fewer subdomains. Thereby, the local solver shows good strong scalability as indicated in the
MBFs column. Interestingly, the Msc-Gain is around $5$ irrespective of the problem
for the first coarse scale decomposition. However, as the burden of the interface system
relatively grows, as is the case for the second decomposition, the Msc-Gain behaves
differently depending on the problem, with values that oscillates between $1.5$ up to $5.7$.
These results suggest the use of configurations with a limited number of interface
unknowns as well as the necessity for better interface linear solvers.}
 
\begin{table}[]
	\centering
	\begin{tabular}{lcrr}
		\toprule
		& coarse decomp            & $60\times10\times50$ & $60\times20\times50$       \\ 
	
		& $N_{\Gamma}$                 & 171.8K           & 349.6K                \\ 
	
		& \#subdomains                 & 30K             & 60K                 \\
		\midrule 
		& fine grid                     & $N_{{\Omega}^\ell}$ & $N_{{\Omega}^\ell}$    \\ 
		$P_1^{ssc}$ & $360\times660\times300$ ({\bf 71M})  & 2.4K                & 1.2K     \\
		$P_2^{ssc}$ & $720\times1320\times600$ ({\bf 570M}) & 19.0K               & 9.5K   \\
		$P_3^{ssc}$ & $660\times2420\times600$ ({\bf 958M}) & 32.0K               & 16.0K  \\ \bottomrule
	\end{tabular}
\caption{Strong scaling setup for three grid configurations and two coarse decompositions.}\label{strong_tab}
\end{table}

\begin{table}[h!!]
  \centering
      {\color{black}
	\begin{tabular}{ccccccccc}
		\toprule
		\multicolumn{1}{l}{}             & \multicolumn{1}{l}{} & Int. time     & \multicolumn{3}{c}{1.7}                  & \multicolumn{3}{c}{7.2}                  \\ \cline{2-9} 
		\multicolumn{1}{l}{}             & \#cores              & Fine    & MRCM    & Msc-Gain    & MBFs    & MRCM    & Msc-Gain    & MBFs     \\ \cline{2-9} 
		\multicolumn{1}{c}{$P_1^{ssc}$} &                      &               &         &          &         &           &         &               \\ \cline{1-1}
		& 150                    & 20.7          & 5.8     & 3.6      & 4.0          & 11.6    & 1.8      & 4.0          \\
		& 300                   & 17.5          & 3.9     & 4.5      & 2.2          & 9.6     & 1.8      & 2.2          \\
		& 600                   & 14.4          & 2.8     & 5.1      & 1.2            & 8.5     & 1.7      & 1.1         \\
		& 1000                   & 12.7          & 2.5     & 5.1      & 0.7          & 8.1     & 1.6      & 0.7             \\ \cline{2-9} 
		\multicolumn{1}{c}{$P_2^{ssc}$} &                      &               &         &          &         &           &         &               \\ \cline{1-1}
		& 150                    & 253.5         & 41.9    & 6.0      & 38.9           & 46.0    & 5.5      & 35.8          \\
		& 300                   & 110.9         & 23.7    & 4.7      & 21.4           & 27.9    & 4.0      & 19.2     \\
		& 600                   & 58.1          & 13.0    & 4.5      & 10.9           & 17.9    & 3.2      & 9.9           \\
		& 1000                   & 41.8          & 8.7     & 4.8      & 6.7            & 13.9    & 3.0      & 6.2       \\ \cline{2-9} 
		\multicolumn{1}{c}{$P_3^{ssc}$} &                      &               &         &          &         &           &         &            \\ \cline{1-1}
		& 150                    & -             & 63.1    & -        & 58.7         & 70.1    & -        & 60.6      \\
		& 300                   & 230.8         & 36.1    & 6.4      & 33.1        & 40.8    & 5.7      & 32.7        \\
		& 600                   & 97.8          & 19.7    & 5.0      & 17.0         & 24.5    & 4.0      & 16.8   \\
		& 1000                   & 68.2          & 12.9    & 5.3      & 10.6         & 18.1    & 3.8      & 10.4       \\ \bottomrule
	\end{tabular}}
\caption{Strong scaling assessment for the configurations defined in Table \ref{strong_tab}.}\label{sts}
\end{table}

\subsubsection{Computational time dependence on the permeability field contrast}

{\color{black}
As a further verification of the relative advantage of our multiscale method based on direct solvers
with respect to a fine grid solver based on iterative methods we investigate the effect
on the computational performance of the permeability contrast, which is a critical parameter
in reservoir engineering. To that end, let us take a parameter $\theta>0$ and define the permeability
tensor $\mathbf{K}_\theta$ by modifying each component of $\mathbf{K}$ according to
$$\mathbf{K}^{ij}_\theta\left(\mathbf{x}\right):=\left(\mathbf{K}^{ij}\left(\mathbf{x}\right)\right)^\theta~.$$
This way, the permeability contrast ($c=\mathbf{K}_{\max}/ \mathbf{K}_{\min}$) is modified as $c\rightarrow c^\theta$.

The modified SPE10 field is then projected on a fine grid consisting of $90$M cells as previously explained.
Computational times are reported for different values of $\theta$ in Table \ref{tab:time-amg-contrast}.
In order to simplify the analysis, the computational times correspond to fixing the
number of GMRES iterations to 12 (by setting \verb|-ksp_max_it 12| and \verb|-ksp_rtol 1e-20|).
This was sufficient for the relative residual to reach a predefined tolerance of $10^{-8}$ in all runs.
Interestingly, for low values of the contrast, the computational time spent for the iterative solver preconditioned
with BoomerAMG decreases as $c$ is increased. However, for high values of $c$. Similar to the SPE10 original contrast
of $\sim 10^7$ and beyond, the computational time is nearly constant.
This is consistent with \cite{Ruge1987}, the seminal work that introduced the algebraic multigrid on which BoomerAMG
is based. According to \cite{Ruge1987} the computational complexity relies on two parameters: on one hand,
$\sigma^\Omega$, denoting the ratio of the total number of points on all grids to that on the fine grid
that is referred to as the grid complexity and on the other hand, $\sigma^A$, denoting the ratio of the total number
of nonzeros entries in all the matrices to that in the fine grid matrix, which is called the operator complexity.
In Table \ref{tab:time-amg-contrast} it is observed that $\sigma^\Omega$ is relatively constant for the tested
configurations, while the variation of $\sigma^A$ is in direct relation to that of the computational time.
As a result, the Msc-Gain varies from $8.2$ for low values of $\theta$ up to something close to
$3.7$ for high values of $\theta$.}

\begin{table}[h]
	\begin{center}
		\begin{tabular}{cccccc}
			\toprule
			$\theta$&	Contrast	& $\sigma^\Omega$ & $\sigma^A$ & time & Msc-Gain \\ 
			\midrule
			0.30 &	$1.9\times 10^2$  & 2.13 & 13.6 & 115.3 & 8.2 \\
			0.48 &	$4.4\times 10^3$  & 2.03 & 9.87 & 75.9 & 5.4 \\
			0.65 &	$8.7\times 10^4$  & 1.97 & 6.57 & 54.1 & 3.9 \\
			0.83 &	$2.0\times 10^6$  &  1.97 & 5.93 & 51.5 & 3.7 \\
			1.0 &	$\mathbf{4.0\times 10^7}$  &  1.98 & 6.05 & 50.4 & 3.6 \\
			1.17 &	$7.8\times 10^8$  & 1.99 & 6.11 & 50.4 & 3.6 \\
			1.35 &	$1.8\times 10^{10}$  & 2.00 & 6.16 & 51.9 & 3.7 \\
			\bottomrule
		\end{tabular}
		\caption{Dependence between the computational time of GMRES preconditioned with BoomerAMG when varying the contrast of the SPE10 benchmark.}\label{tab:time-amg-contrast}
	\end{center}
\end{table}

\subsubsection{Choice of multiscale method}

\begin{figure}[h!!]
	\includegraphics[width=\textwidth]{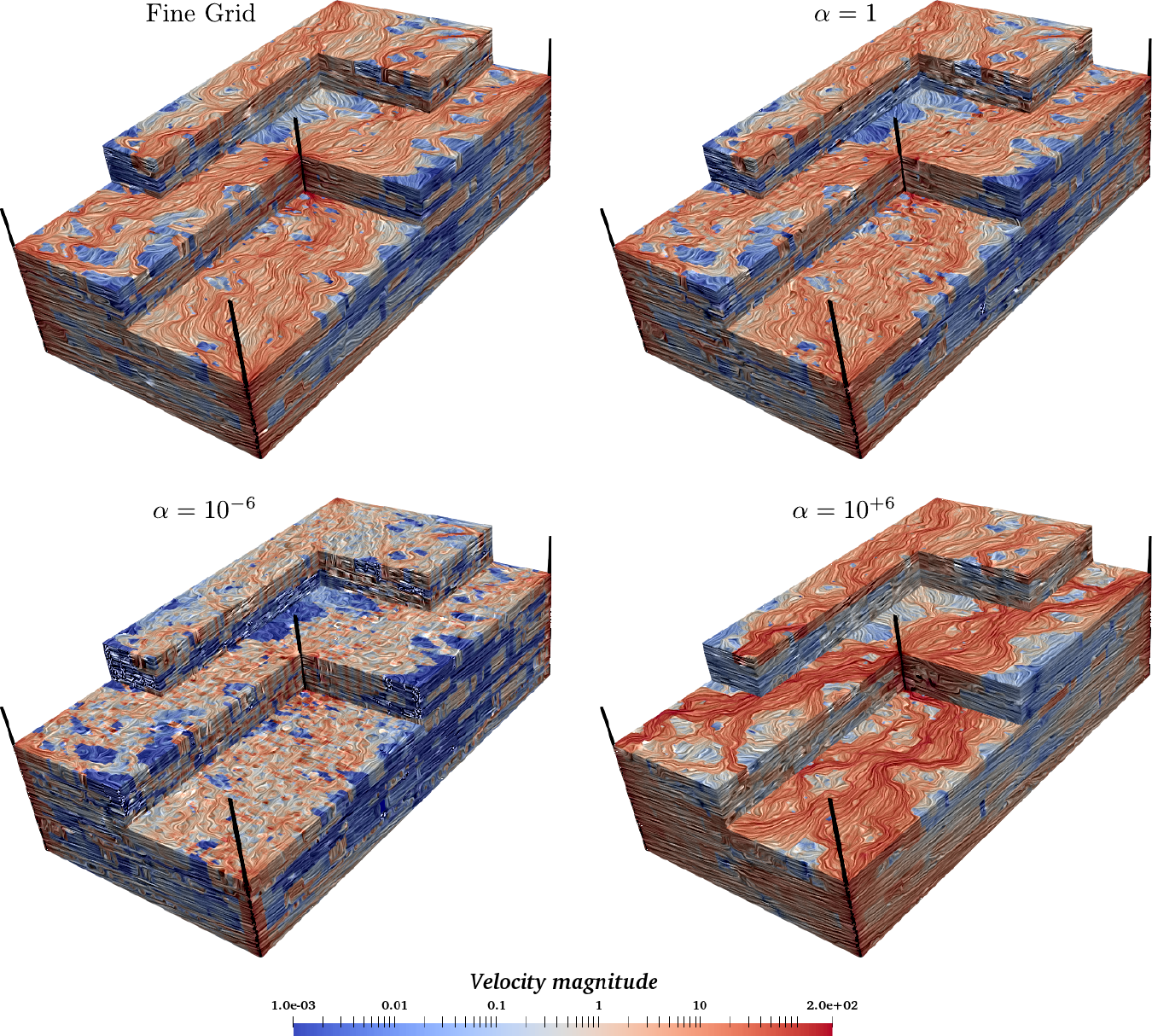}
	\caption{Fine grid solution and the multiscale solution for different choices of the $\alpha$ parameter .}\label{fig:comparison-alpha-perspective}
\end{figure}

\begin{figure}[h]
	\centering
	\includegraphics[width=0.8\textwidth]{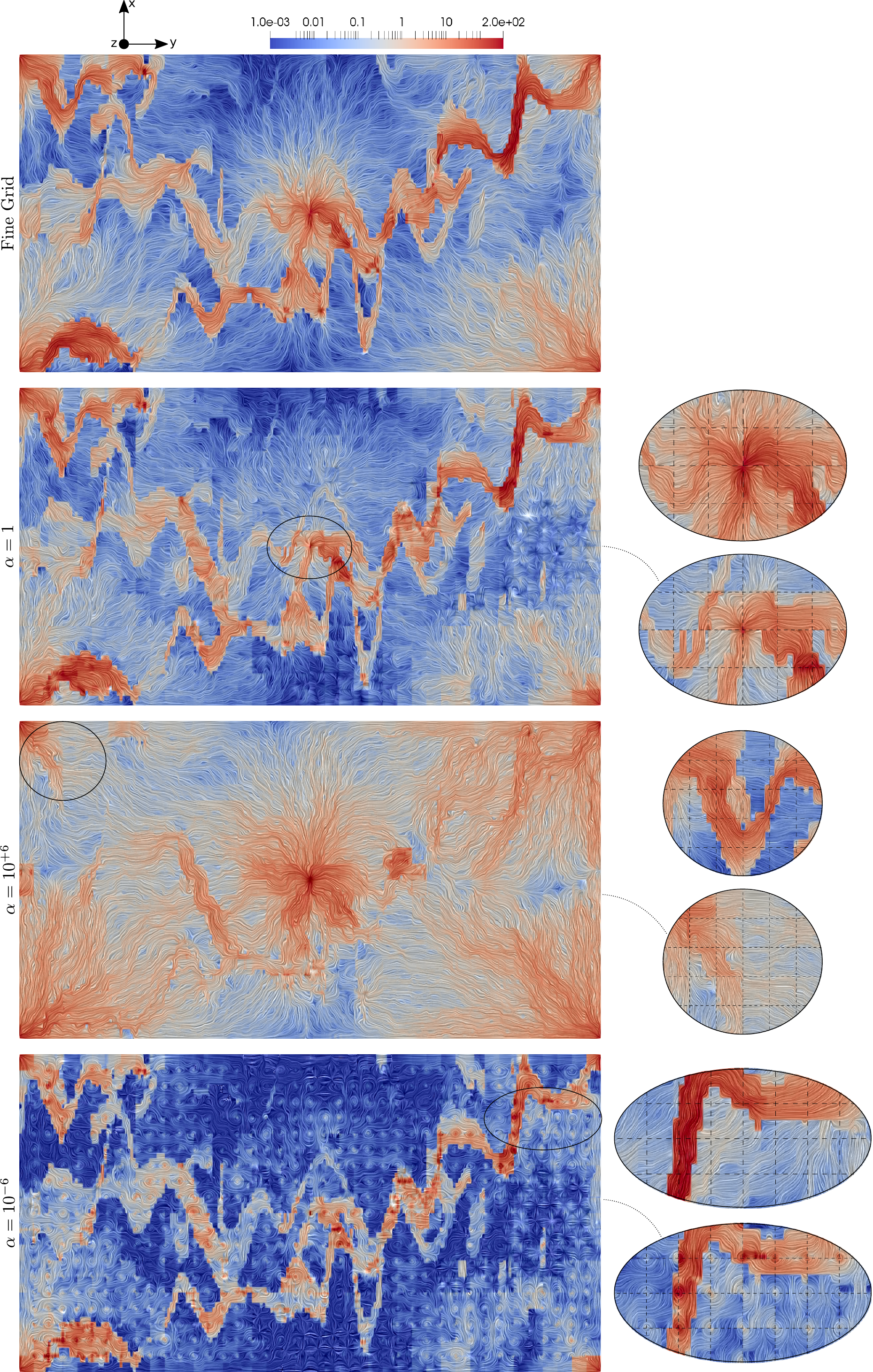}
	\caption{Fine grid solution and the multiscale solution for different choices of the $\alpha$ parameter for a cut in layer 36 in $x_3$ plane of the 3D solution.}\label{fig:comparison-alpha-Z2}
\end{figure}

\begin{figure}[h]
	\centering
	\includegraphics[width=0.6\textwidth]{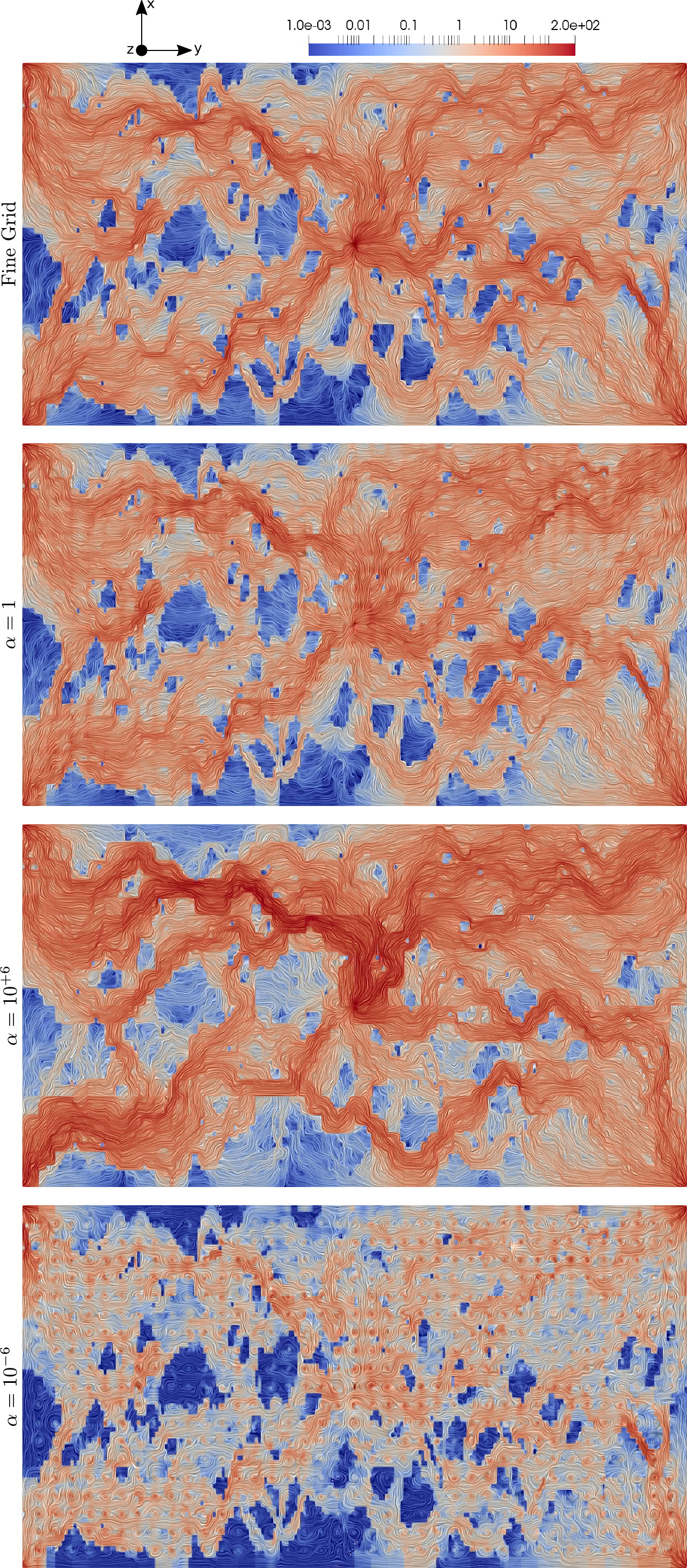}
	\caption{Fine grid solution and the multiscale solution for different choices of the $\alpha$ parameter for a cut in layer 85 in $x_3$ plane of the 3D solution.}\label{fig:comparison-alpha-Z3}
\end{figure}

In this article we have chosen the MRCM as a prototypical formulation from which other multiscale methods can be obtained
by varying the interface Robin condition parameter $\alpha$ in (\ref{Robin_par}). If we consider the limit $\alpha \to 0_{+}$,
the multiscale solution satisfies pointwise pressure continuity at the skeleton $\Gamma$ of the domain decomposition.
On the other hand, as $\alpha \to \infty$, the multiscale solution satisfies pointwise velocity continuity at $\Gamma$.
The numerical experiments presented so far were obtained for $\alpha = 1$, a choice
based on prior numerical experience, for which neither the pressure nor the velocity are continuous
on $\Gamma$. The velocity field produced by the multiscale method for such case is depicted in the top part
Fig. \ref{fig:comparison-alpha-perspective} jointly with the fine grid solution. 
The multiscale solution corresponding to the extreme values $\alpha=10^{-6}$ and $10^{+6}$
are displayed in the bottom part.
The multiscale solution corresponding to $\alpha = 1$ captures remarkably well the underlying features
of the permeability field showing great similarity to the fine grid solution. Numerical results
for $\alpha=10^{-6}$ and $10^{+6}$ are significantly far away from those. This is better noticed in
in Figs. \ref{fig:comparison-alpha-Z2} and \ref{fig:comparison-alpha-Z3}, in which the 3D solutions are intersected
by planes corresponding to layers 36 and 85 of the SPE10, precisely where interesting highly channelized structures are located.
The reference solution for layer 36 is displayed at the top of Fig. \ref{fig:comparison-alpha-Z2} followed by the multiscale solutions.
The inserts show zoomed regions of interest, one for each choice of $\alpha$, pared to the corresponding area in the reference
solution to ease the comparison.
For the MRCM($10^{+6}$), the zoomed area displays a region having a high permeability channel and a low permeability
background. We notice in such case the method does not detect the presence of the channel and the fluid passes through it.
For the MRCM($10^{-6}$) the zoomed area also points to a channel with high permeability.
Whilst the channelized structure is better represented, unphysical recirculation patterns are observed
near the corners of the domain decomposition. Similar patterns are also observed in Fig. \ref{fig:comparison-alpha-Z3}
that corresponds to layer 85. More sophisticated choices of the algorithmic parameter as those proposed
in \cite{rocha2020} can be adopted with the potential to deliver more accurate results at the same computational cost.

\begin{figure}[h]
	\centering
	\includegraphics[width=0.95\textwidth]{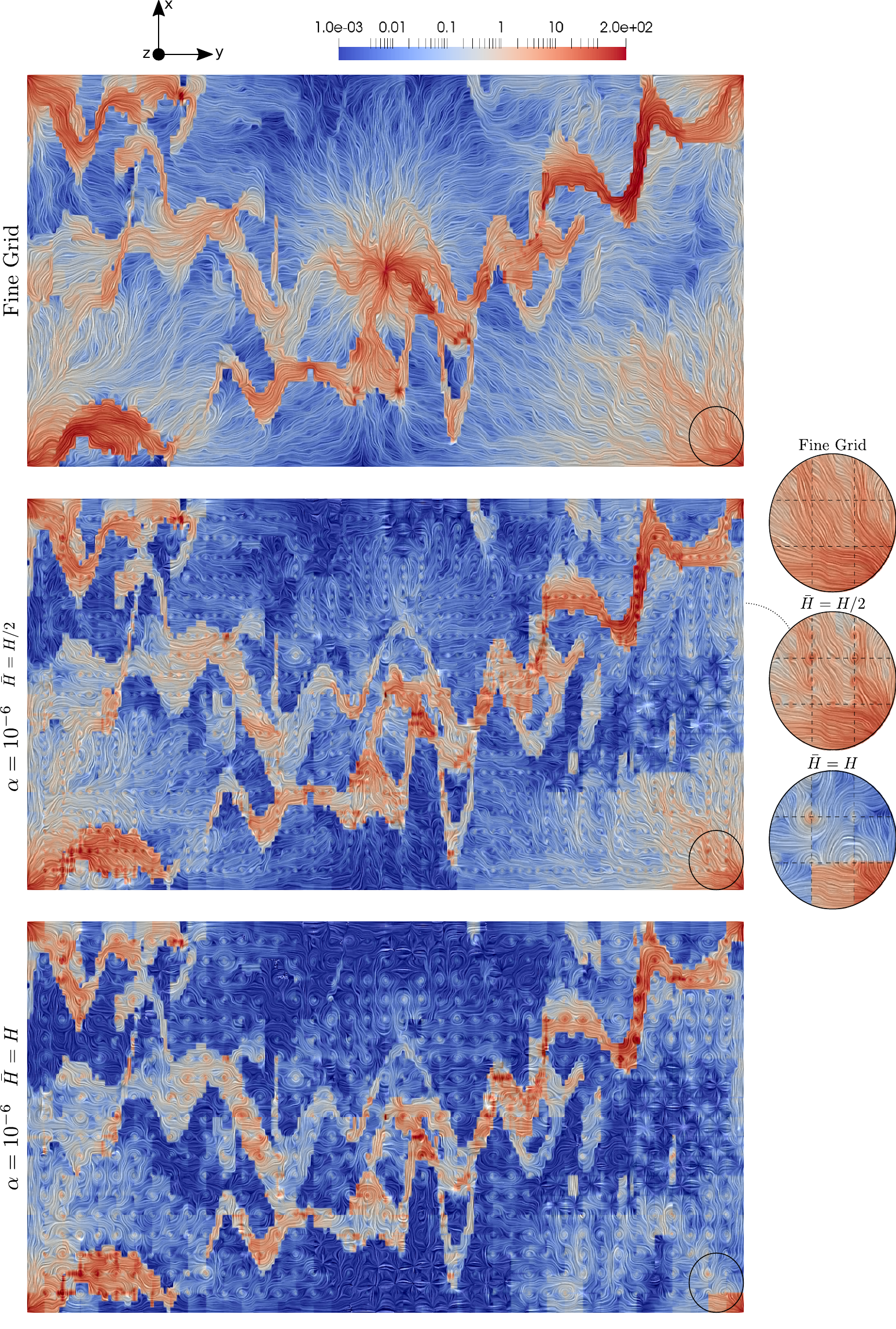}
	\caption{Fine grid solution and the multiscale solution for $\alpha=10^{-6}$ and two values of $\bar{H}$ for a cut in layer 36 along the $x_3$ plane of the 3D solution.}\label{fig:Z2-refhbar2}
\end{figure}

{\color{black}
Also, accuracy can be improved by making smarter choices of the interface spaces
defined over the skeleton $\Gamma$. We have limited so far to the piecewise constant case.
In our implementation, nonetheless, we can easily perform an $H$-refinement over $\Gamma$
as explained in \ref{sec:MRCMexplained}. To illustrate this, Fig. \ref{fig:Z2-refhbar2} shows the solution
for the MRCM($10^{-6}$) obtained by halving the interface elements (i.e., $\bar{H} = H/2$).
Although the dimension of the interface coupling system multiplies by 4, the
computation of multiscale basis function comes almost at no cost, since it only
involves the solution of a linear system for a greater number of right hand sides.
One can notice an overall improvement of the solution especially in regions of high-speed
magnitude as occurs close to wells (see the zoomed areas in the inserts).

Further improvements could make use of informed functions spaces or spaces based on physics
(see \cite{guiraldello2019,rochaphd}) that can better capture
the underlying variations of the rock formation, with possibly a reduced number of
degrees of freedom.}

\subsubsection{Two-phase flow simulations} 

\label{sec:hpc-biphasic}
\input{tex/hpc_biphasic.tex}

\section{Conclusions and outlook}
\label{sec:conclusions}

We have presented an HPC implementation for multiscale mixed methods based on non-overlapping domain decomposition
techniques, especially tailored for the simulation of flow in porous media.
We considered problems ranging from a
several million up to a billion fine grid cells. The accuracy and computational performance of the scheme
as well as its scalability properties have been assessed in problems involving anysotropic highly heterogeneous
permeability fields based on the SPE10 benchmark. For testing of different multiscale options the general framework
offered by the MRCM method have been adopted.

The main conclusions that emerge from exhaustive numerical experimentation are:

\begin{itemize}

\item The best computational times are obtained for the multiscale method
  by using direct solvers for both the local problems and the interface 
  system. In the former case, this comes from the fact that several right hand sides
  must be solved for each subdomain. In the latter case, this comes from
  the conditioning of the coupling system. 

\item Assessment of the direct solver used in this work (MUMPS) shows that
  decreasing the number of local unknowns below a certain treshold $N_{\Omega^\ell}\sim 40$K,
  yielding more subdomains per processor, will eventually lead to a plateau in
  computing time of the MBFs, since the behavior changes from quadratic
  to linear in $N_{\Omega^\ell}$ at this point.
  
\item In general, the multiscale method exhibits better weak scalability and computational performance for
  the problem sizes considered in this work. This depends on the number of local 
  and interface unknowns that can be chosen as the domain decomposition is specified.
  This good behavior is limited by the size of the interface system. For values of $N_{\Gamma} \gtrsim 200$K
  a decline in scalability properties of the scheme is observed. 

\item The strong scalability assessment reveals that if the number of processors available is
  a limiting factor, the multiscale mixed method is an excellent asset, either because the problem may
  not fit in memory for the fine grid solver or due to the significant gain in computing time obtained,
  that can be as high as $6$.
  
\item One advantage of the strategy adopted is that it allows to estimate a priori the computing
  time of the solver irrespective of the problem and permeability contrast. 

\item Accuracy of the multiscale solver has been also assessed in a two phase flow scenario
  in a problem involving $406$M unknowns. Results for the production curves and water cuts are
  promissory yielding errors smaller that $6\%$ with respect to the fine grid solution.
  
\end{itemize}

Results presented so far show the potential of Multiscale Mixed Methods to solve
large scale porous media problems.

\section*{Acknowledgments}

The authors gratefully acknowledge the financial support received from the Brazilian oil company Petrobras grant 2015/00400-4,
from the S\~ao Paulo Research Foundation FAPESP/CEPID/CeMEAI grant 2013/07375-0, and from Brazilian National Council for Scientific and Technological Development
CNPq grants 305599/2017-8 and 310990/2019-0. This research was carried out using computational resources from the Cluster Euler,
Centre for Mathematical Sciences Applied to Industry (CeMEAI).

\bibliographystyle{unsrt}

\end{document}

%% file: tex/hpc_math_num_models.tex
Consider a domain $\Omega\subset \mathbb{R}^3$ and two-phase flow of water and oil with respective saturations fields denoted by
$S_\text{w}$ and $S_\text{o}$, that depend both on position $\mathbf{x}$ and time $t$. Capillary pressure and gravity effects are neglected.
Assuming saturated porous media we have the restriction
\begin{equation}
S_\text{w}+S_\text{o}=1.
\end{equation}
The evolution of $S_\text{w}$ is modeled by the Buckley-Leverett equation
\begin{equation}
\parder{\Sw}{t}+\nabla\cdot \left( \varphi(\Sw) \mathbf{u} \right)=0\,,
\label{eq:sec2:Buckley-Leverett}
\end{equation}
where $\varphi$ is the so called fractional flow function. By rescaling the time variable the porosity of the medium, 
which is assumed constant in this work, is removed from eq. (\ref{eq:sec2:Buckley-Leverett}).

Velocity $\mathbf{u}$ and pressure field $p$ satisfy the Darcy's law and the continuity equation with prescribed boundary conditions
\begin{align}
\mathbf{u} & =  -\lambda(\Sw) \mathbf{K}\,\nabla p && \text{in }\Omega \label{eq:sec2:darcy}\\	
\nabla \cdot \mathbf{u} & = f && \text{in }\Omega
\label{eq:sec2:continuity}
\\	
p & =  g_p && \text{on }\partial\Omega_p
\label{eq:sec2:bcp}
\\	
\mathbf{u}\cdot \hat{\mathbf{n}} & =  g_u && \text{on }\partial\Omega_u
\label{eq:sec2:bcu}
\end{align}
where $\mathbf{K}$ is the permeability tensor assumed to be diagonal and anisotropic, i.e., $\mathbf{K}=\text{diag}\,(K^{11},K^{22},K^{33})$, $\lambda$ is the total mobility function, $g_p$, $g_u$ and $f$ are given {\it forcing} data. The functions $\varphi$ and $\lambda$ read
$$\lambda(\Sw) = \frac{\Sw^{2}}{\mu_\text{w}} + \frac{(1-\Sw)^{2}}{\mu_\text{o}} \quad \text{and}\quad \varphi(\Sw) = \frac{(\mu_\text{o}/\mu_\text{w})\Sw^{2}}{(\mu_\text{o}/\mu_\text{w})\Sw^{2} + (1-\Sw)^{2}}\,.$$
in which $\mu_\text{w}$ and $\mu_\text{o}$ stands for the viscosity of the water and oil phases, respectively.

\subsection{Domain decomposition and Darcy solvers} \label{sec:elliptic-solvers}
The domain is decomposed into non-overlapping subdomains $\Omega^{\ell}$ of typical size $H$,
such that $\Omega = {\bigcup}_{\ell}\, \Omega^{\ell}$. The skeleton of the domain decomposition
is denoted by $\Gamma=\bigcup_{\ell}\,\Gamma_{\ell,\ell'}$, with $\Gamma_{\ell,\ell'}=\bar{\Omega}^{\ell}\cap \bar{\Omega}^{\ell'}$.
Let us denote $\hat{\bf n}^{\ell}$ a unit vector normal to $\Omega^{\ell}$,
define $\hat{\mathbf{n}}=\hat{\mathbf{n}}^{\min\{\ell,\ell'\}}$ and use $+$ and $-$ superscripts
to refer to the two-sided limits approaching $\Gamma$.

Two different Darcy solvers are implemented to approximate the solution to \eqref{eq:sec2:darcy}-\eqref{eq:sec2:bcu}, namely,
a fine grid global method and the MRCM.
The main difference between these two solvers is that continuity of the unknown fields
$\mathbf{u}$ and $p$ on the skeleton of the domain decomposition is relaxed in the MRCM by imposing
\begin{align}
\int_\Gamma \left ( {\bf u}^+ -
{\bf u}^- \right )\cdot \hat{\bf n}~v_H \, d{\Gamma}&=0~,
\qquad \forall\,v_H\,\in\,\mathcal{P}_H~, \label{eq:sec2:contflux} \\
\int_\Gamma (p^+-p^-)\,w_H \, d{\Gamma}&=0~,\qquad \forall\,w_H\,\in\,\mathcal{U}_H~, \label{eq:sec2:contpres}
\end{align}
where $\mathcal{P}_H$ and $\mathcal{U}_H$, are suitably defined spaces over the skeleton $\Gamma$ such that
mass conservation is satisfied at a coarse scale $\bar{H}$. By setting $\mathcal{P}_H =\mathcal{P}_h$ and
$\mathcal{U}_H =\mathcal{U}_h$ the solutions of both methods coincide.

A key ingredient for both solvers is a Finite Volume discretization method that, for completeness, is recalled next.
 For any subdomain $\omega\subset \Omega$, we consider a hexahedral mesh to discretize
\eqref{eq:sec2:continuity} by considering \eqref{eq:sec2:darcy} along two-point formulas for the fluxes
as in \cite{guiraldello2018multiscale}.
Consider computational cells of sizes $h_{x_1}$, $h_{x_2}$, $h_{x_3}$ in each direction,
$h=\max\,\{h_{x_1},h_{x_2},h_{x_3}\}$ denotes the maximum of $h_{x_j}$
and the total number of cells is denoted by $N^\omega$. The product of $\lambda$ and $\mathbf{K}$ in Darcy's law \eqref{eq:sec2:darcy}
is referred to as $\mathbf{\hat{K}}$. With this, one ends up with the system of equations
\begin{equation}
\left \{
\begin{array}{r }
-~a_{i-1,j,k}\,p_{i-1,j,k}-a_{i,j,k}\,p_{i+1,j,k}-b_{i,j-1,k}\,p_{i,j-1,k}-b_{i,j,k}\,p_{i,j+1,k}\\
~\\
~~~~-~c_{i,j,k-1}\,p_{i,j,k-1}-c_{i,j,k}\,p_{i,j,k+1}+d_{i,j,k}\,p_{i,j,k}=\,f_{i,j,k}\,,
\end{array}
\right. \label{eq:sec2:darcy-discrete}
\end{equation}
where
\begin{equation}
\begin{cases}
a_{i,j,k}=\dfrac{\tilde{K}^{11}_{i+\frac{1}{2},j,k}}{h_{x_1}^2},~~
\,b_{i,j,k}=\dfrac{\tilde{K}^{22}_{i,j+ \frac{1}{2},k}}{h_{x_2}^2},\,~~c_{i,j,k}=\dfrac{\tilde{K}^{33}_{i,j,k+ \frac{1}{2}}}{h_{x_3}^2}\,,\\
~~\\
\vspace*{0.8mm}
d_{i,j,k}
=\dfrac{\tilde{K}^{11}_{i+\frac{1}{2},j,k}}{h_{x_1}^2}
+\dfrac{\tilde{K}^{{11}}_{i-\frac{1}{2},j,k}}{h_{x_1}^2}
+\dfrac{\tilde{K}^{{22}}_{i,j+\frac{1}{2},k}}{h_{x_2}^2}
+\dfrac{\tilde{K}^{22}_{i,j-\frac{1}{2},k}}{h_{x_2}^2}
+\dfrac{\tilde{K}^{{33}}_{i,j,k+\frac{1}{2}}}{h_{x_3}^2}
+\dfrac{\tilde{K}^{33}_{i,j,k-\frac{1}{2}}}{h_{x_3}^2}\,,\\
~~\\
\vspace*{0.8mm}
\tilde{K}^{11}_{i\pm\frac{1}{2},j,k}=\dfrac{2\,\hat{K}^{11}_{i,j,k}\,\hat{K}^{11}_{i\pm1,j,k}}{\hat{K}^{11}_{i,j,k}+\hat{K}^{11}_{i\pm1,j,k}},~~
\tilde{K}^{22}_{i,j\pm\frac{1}{2},k}=\dfrac{2\,\hat{K}^{22}_{i,j,k}\,\hat{K}^{22}_{i,j\pm1,k}}{\hat{K}^{22}_{i,j,k}+\hat{K}^{22}_{i,j\pm1,k}},~~
\tilde{K}^{33}_{i,j,k\pm\frac{1}{2}}=\dfrac{2\,\hat{K}^{33}_{i,j,k}\,\hat{K}^{33}_{i,j,k\pm1}}{\hat{K}^{33}_{i,j,k}+\hat{K}^{33}_{i,j,k\pm1}},
\end{cases}
\label{eq:sec2:aijk}
\end{equation}
and $i,j,k$ the cells indices along the $x_1$, $x_2$ and $x_3$ axis.
Also, in (\ref{eq:sec2:aijk}), $\tilde{K}^{ll}_{i,j,k}$ terms correspond to the so called harmonic mean of the permeability diagonal entries of $\hat{K}$.
To account for the boundary conditions some of these terms for cells
having a face on the boundary $\partial{\Omega}$ must be modified
(see \cite{guiraldello2018multiscale}).

\subsubsection{The global fine grid solver}
Assembling eqs. \eqref{eq:sec2:darcy-discrete}-\eqref{eq:sec2:aijk} for the case $\omega=\Omega$ 
leads to a linear system given by
\begin{equation}
\mathcal{A}^{h}\,\mathbf{p}^{h}=\mathbf{f}^{h}~,
\label{eq:sec2-darcy-discrete-fine-grid}
\end{equation}
where $\mathcal{A}^{h} \in \mathbb{R}^{N^{\Omega}\times N^{\Omega}}$ represents a sparse matrix operator and $\mathbf{p}^h \in \mathbb{R}^{N^{\Omega}}$
is the global vector of cell pressure unknowns whit $N^{\Omega}$ standing for the total number of cells in the global (target) problem.
 Once this discrete pressure field is found the discrete velocity
field $\mathbf{{u}}^h$ can be recovered at the cell faces by means of two-point formulas, e.g.,
\begin{equation}
u_{i+\frac{1}{2},j,k} = - \tilde{K}^{11}_{i+\frac{1}{2},j,k}\dfrac{p_{i+1,j,k} - p_{i,j,k}}{h_{x_1}},
\label{eq:two_point}
\end{equation}
with analogous definitions for other cell fluxes.
Typical sizes of the global systems to be solved in this work varies between a few million up to a few billion cells,
that precludes the use of direct solvers.

\subsubsection{A domain decomposition Multiscale Mixed method} \label{sec:MRCMexplained}
In this section we recall the basic ingredients of the multiscale method adopted in this work, the MRCM,
that has been presented in detail in \cite{guiraldello2018multiscale,guiraldello2019}.
This multiscale strategy localizes the problem \eqref{eq:sec2:darcy}-\eqref{eq:sec2:bcu} into each subdomain $\Omega^{\ell}$. The boundary conditions \eqref{eq:sec2:bcp}-\eqref{eq:sec2:bcu} are imposed at the regions of $\partial \Omega^\ell$ that lay on  $\partial \Omega$. To fully define the local problems, Robin boundary conditions are imposed at the regions of $\partial \Omega^\ell$ that intercept $\Gamma$. One arrives to a multiscale solution in $\Omega$ by coupling the local solutions through the weak compatibility conditions \eqref{eq:sec2:contflux} and \eqref{eq:sec2:contpres}.
The differential formulation of this method reads: Find subdomain solutions ($\mathbf{u}^{\ell},p^{\ell}$), decomposed as 
$
\mathbf{u}^{\ell}=\mathbf{\hat{u}}^{\ell}+\mathbf{\bar{u}}^{\ell},~{p}^{\ell}={\hat{p}}^{\ell}+{\bar{p}}^{\ell}
$
and interface fields $(U_H,P_H)$ satisfying:

\begin{equation}
\left.
\begin{array}{r l l l l}
{\bf \bar{u}}^{\ell} & = & -\, \mathbf{K}\,\nabla \bar{p}^{\ell} & \mbox{in}~\Omega^{\ell} \\
\nabla \cdot {\bf \bar{u}}^{\ell} & = & f & \mbox{in}~\Omega^{\ell}\\
\bar{p}^{\ell} & = & g_p & \mbox{on}~\partial{\Omega^{\ell}} \cap \partial\Omega_p\\
{\bf \bar{u}}^{\ell} \cdot \hat{\bf n}^{\ell} & = & g_u & \mbox{on}~\partial{\Omega^{\ell}} \cap \partial\Omega_u\\
-\beta^{\ell} {\bf \bar{u}}^{\ell} \cdot \hat{\bf n}^{\ell} + \bar{p}^{\ell} & = & 0 & \mbox{on}~ \partial{\Omega^{\ell}} \cap \Gamma\\
\end{array}
\right\} \label{eq:sec2:localomegapar}
\end{equation}

\begin{equation}
\left .
\begin{array}{r l l l l}
{\bf \hat{u}}^{\ell} & = & -\, \mathbf{K}\,\nabla p^{\ell} & \mbox{in}~\Omega^{\ell} \\
\nabla \cdot {\bf \hat{u}}^{\ell} & = & 0 & \mbox{in}~\Omega^{\ell}\\
\hat{p}^{\ell} & = & 0 & \mbox{on}~\partial{\Omega^{\ell}} \cap \partial\Omega_p\\
{\bf \hat{u}}^{\ell} \cdot \hat{\bf n}^{\ell} & = & 0 & \mbox{on}~\partial{\Omega^{\ell}} \cap \partial\Omega_u\\
-\beta^{\ell} {\bf \hat{u}}^{\ell} \cdot \hat{\bf n}^{\ell} + \hat{p}^{\ell} & = & -\beta^{\ell} U_H \, \hat{\bf n} \cdot \hat{\bf n}^{\ell} + P_H & \mbox{on}~ \partial{\Omega^{\ell}} \cap \Gamma\\
\end{array}
\right\} \label{eq:sec2:localomegaint}
\end{equation}
\begin{equation}
\displaystyle \sum _{\ell}  { \int _{\partial\Omega^{\ell}\cap\Gamma} \left ( {\bf u}^{\ell}\cdot \hat{\mathbf{n}}^{\ell} \right ) \, \phi \, d{\Gamma}} = 0\qquad \forall\,\phi\,\in\,\mathcal{P}_H \label{eq:sec2:weakcont_f}
\end{equation}
\begin{equation}
\displaystyle \sum _{\ell} \int _{\partial\Omega^\ell\cap\Gamma}  \beta^{\ell} \left ({\bf u}^\ell\cdot \hat{\mathbf{n}}^{\ell} -  U_H \, \hat{\mathbf{n}}^{\ell} \cdot \hat{\mathbf{n}} \right ) \hat{\mathbf{n}}^{\ell}\cdot\hat{\mathbf{n}}~\phi \,d{\Gamma}= 0\qquad \forall\,\phi\,\in\,\mathcal{U}_H \label{eq:sec2:weakcont_p}
\end{equation}
where $\beta^\ell$ are the Robin condition parameters.
The differential problems \eqref{eq:sec2:localomegapar} and \eqref{eq:sec2:localomegaint} are solved independently.
In addition an interface linear problem that results from \eqref{eq:sec2:weakcont_f}-\eqref{eq:sec2:weakcont_p}
must be solved. The interface problem can be written by first
introducing a finite dimensional space $\mathcal{V}_H = \mbox{span}\{\phi_1,\dots,\phi_{N_\mathcal{V}}\}$ defined on the
skeleton of the domain decomposition such that the interface fields read
\begin{equation}
	P_H=\sum_{k=1}^{N_\mathcal{V}}\pi^p_k\,\phi_k~,\qquad U_H=\sum_{k=1}^{N_\mathcal{V}}\pi^\mathbf{u}_k\,\phi_k~.
	\label{eq:P_H_U_H}
\end{equation}

\begin{remark} Although the MRCM allows for the use different spaces for $P_H$ and $U_H$ \cite{guiraldello2018multiscale},
in this work, for simplicity, we restrict ourselves to the case where both spaces are the same.
Several choices of interface spaces have been proposed in \cite{guiraldello2018multiscale,guiraldello2019,rochaphd}.
\end{remark}

By denoting $X=(\pi^p_1,\ldots,\pi^p_{N_\mathcal{V}},\pi^\mathbf{{u}}_1,\ldots,\pi^\mathbf{{u}}_{N_\mathcal{V}})^\intercal$
allows us to write the conditions \eqref{eq:sec2:weakcont_f}-\eqref{eq:sec2:weakcont_p} as the linear system 
\begin{equation}
\mathcal{A}^{\text{int}}\, X = \mathbf{b}^{\text{int}} \rightarrow
\begin{bmatrix}
~A^{\tiny{PP}} &A^{\tiny{PU}}~ \\
& \\
~A^{\tiny{UP}} & A^{\tiny{UU}}~
\end{bmatrix} 
\begin{bmatrix}
\pi^p \\
\\
\pi^\mathbf{{u}}
\end{bmatrix} = 
\begin{bmatrix}
{b}^{\tiny{P}} \\
\\
{b}^{\tiny{U}} 
\end{bmatrix},\label{eq:interface_problem}
\end{equation}
where
\begin{equation}
	A^{PP}_{ij}=\sum_{\underset{\ell<\ell'}{\ell,\ell'}}\int_{\Gamma_{\ell,\ell'}}\{\mathbf{\hat{u}}_h(\phi_j,0)\}\,\phi_i\,d\Gamma,~ 	A^{PU}_{ij}=\sum_{\underset{\ell<\ell'}{\ell,\ell'}}\int_{\Gamma_{\ell,\ell'}}\{\mathbf{\hat{u}}_h(0,\phi_j)\}\,\phi_i\,d\Gamma,~ b^P_{i}=\sum_{\underset{\ell<\ell'}{\ell,\ell'}}\int_{\Gamma_{\ell,\ell'}}\{\mathbf{\bar{u}}_h\}\,\phi_i\,d\Gamma. \label{eq:sec2-int-APP}
\end{equation}
In (\ref{eq:sec2-int-APP}) $\widehat{\bf u}_h^\ell(\phi_j,0)$ and $\widehat{\bf u}_h^{\ell}(0,\phi_j)$ is the velocity part of the so called 
Multiscale Basis Functions (MBFs) (see \cite{guiraldello2018multiscale}) obtained as solutions to subdomain problem
\eqref{eq:sec2:localomegaint} by taking $(P_H,U_H) = (\phi_j,0)$ and $(P_H,U_H) = (0,\phi_j)$, respectively.
Similarly, $\hat{p}_h^\ell\left(\phi_j,0\right)$ and $\hat{p}_h^\ell\left(0,\phi_j\right)$
denotes the corresponding pressure part of these local problems.
Also, $\mathbf{\bar{u}}_h^\ell$ are the particular solutions of the local problems \eqref{eq:sec2:localomegapar}.
In \eqref{eq:sec2-int-APP}, $\{\cdot\}$ stands for the jump operator of a vector
field $\mathbf{v}$ on $\Gamma$, i.e., $\{\mathbf{v}\}=(\mathbf{v}^\ell-\mathbf{v}^{\ell'})\cdot\hat{\bf n}$
at each $\Gamma_{\ell,\ell'}$. For the second row block of \eqref{eq:interface_problem} one has:
\begin{equation}
A^{UP}_{ij}=\sum_{\underset{\ell<\ell'}{\ell,\ell'}}\int_{\Gamma_{\ell,\ell'}}\left(\beta_i\mathbf{\hat{u}}_h(\phi_j,0)+\beta_j\mathbf{\hat{u}}_h(\phi_j,0)\right)\cdot{\hat{\bf n}}\,\phi_i\,d\Gamma,~~
\end{equation}\label{eq:sec2-int-AUP}
 \begin{equation}	A^{UU}_{ij}=\sum_{\underset{\ell<\ell'}{\ell,\ell'}}\int_{\Gamma_{\ell,\ell'}}\left[\left(\beta_\ell\mathbf{\hat{u}}_h^\ell(0,\phi_j)+\beta_{\ell'}\mathbf{\hat{u}^{\ell'}}_h(0,\phi_j)\right)\cdot{\hat{\bf n}}-\left(\beta_\ell+\beta_{\ell'}\right)\phi_j\right]\phi_i\, d\Gamma \label{eq:sec2-int-AUU}
\end{equation}
\begin{equation}
	b^U_{i}=-\sum_{\underset{\ell<\ell'}{\ell,\ell'}}\int_{\Gamma_{\ell,\ell'}}\left[\beta_\ell\mathbf{\bar{u}}_h^\ell+\beta_{\ell'}\mathbf{\bar{u}^{\ell'}}_h\right]\cdot{\hat{\bf n}}\,\phi_i\, d\Gamma\label{eq:sec2-int-bU}
\end{equation}
Once $P_H$ and $U_H$ are obtained by solving \eqref{eq:interface_problem}, the multiscale solution is locally recovered by means of \eqref{eq:sec2:localomegaint}.
Notice that one has two options to perform this last step. The first one is just computing the combinations
       	\begin{equation}
        	p^\ell=\sum_{k=1}^{N_\mathcal{V}}\left [\pi^p_k\,{\hat{p}}_h^\ell\left(\phi_k,0\right)  + \pi^{\bf u}_k\,{\hat{p}}_h^\ell\left(0,\phi_k\right) \right ] + {\bar{p}}_h,\qquad
                \mathbf{u}^\ell=\sum_{k=1}^{N_\mathcal{V}}\left [ \pi_k^p\,\mathbf{\hat{u}}_h^\ell\left(\phi_k,0\right) +  \pi_k^\mathbf{{u}}\,\mathbf{\hat{u}}_h^\ell\left(0,\phi_k\right) \right ] + \mathbf{\bar{u}}_h,
        \end{equation}
if the information is still available in memory. Otherwise, we can use \eqref{eq:P_H_U_H} and explicitly
solve problems \eqref{eq:sec2:localomegapar} and \eqref{eq:sec2:localomegaint}.

\begin{figure}[h!]
	\centering
	\def\svgwidth{\linewidth}
	\includegraphics[width=0.8\linewidth]{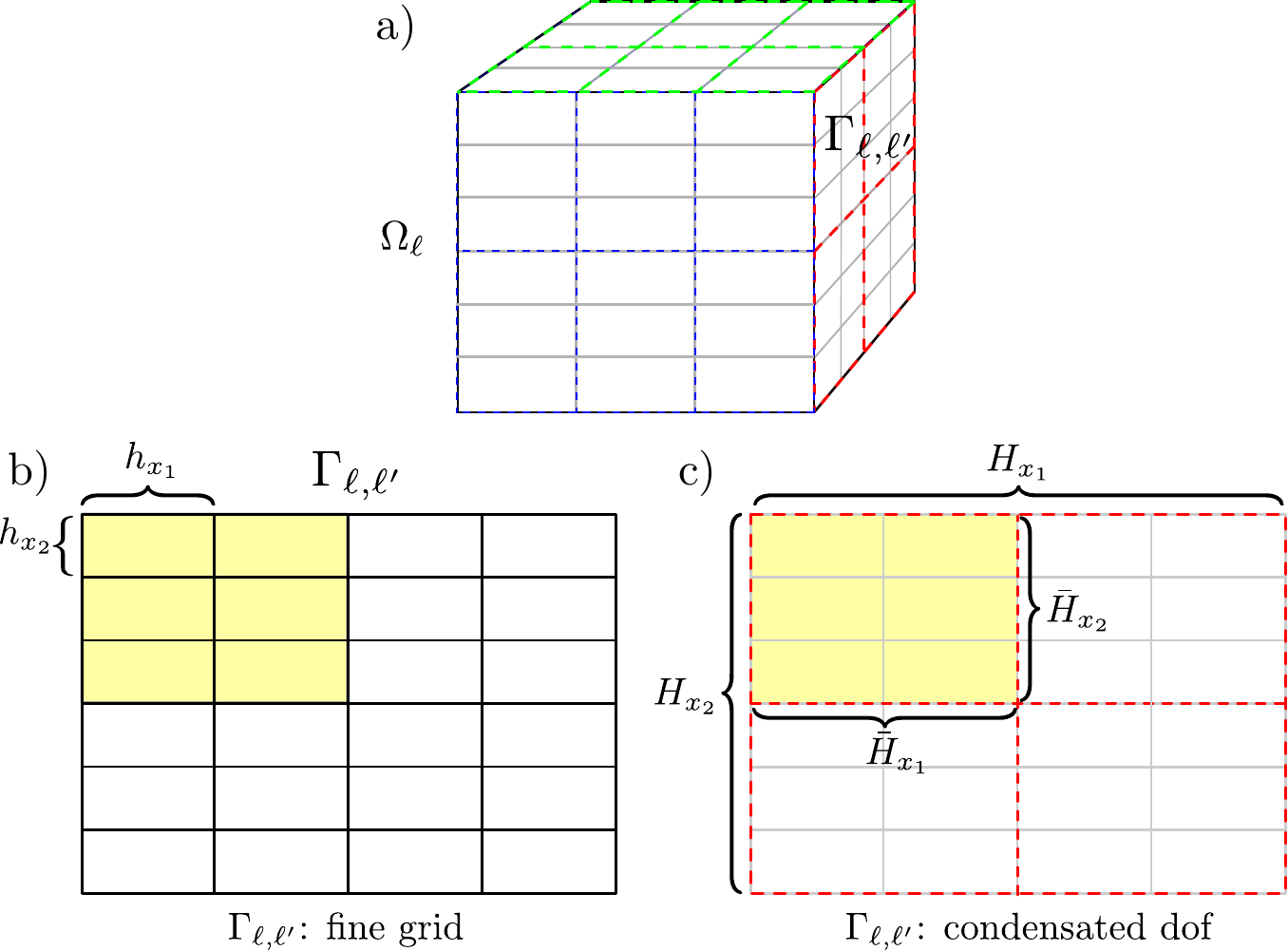}
	\caption{Example showing the condensate degrees of freedom at the interfaces between subdomains.
        In this particular case  the number of degrees of freedom is decreased by a factor of
        ${\frac{\bar{H}_{x_1}}{h_{x_1}}}\frac{\bar{H}_{x_2}}{h_{x_2}}=2\cdot 3=6$.}
	\label{fig:sec2:decomposition_h_H}
\end{figure}

We still need to define the space $\mathcal{V}_H$. To do so, let $H_{x_i}$ ($i=1,2,3$) be the length of
the subdomains $\Omega^\ell$ (integer multiple of the cells length $h_{x_i}$) along some direction and $\bar{H}_{x_i}$ be a multiple of
$h_{x_i}$ such that $h_{x_i}\leq\bar{H}_{x_i}\leq {H}_{x_i}$, $i=1,2,3$. One introduces first the fine grid space $\mathcal{V}_h$ on $\Gamma$.
To that end, consider some arbitrary interface $\Gamma_{\ell,\ell'}$, as shown in Figure \ref{fig:sec2:decomposition_h_H} a) and b),
the intersection of the volumetric cells with $\Gamma_{\ell,\ell'}$ induces a decomposition of $\Gamma_{\ell,\ell'}$ into rectangles of
lengths $h_{x_i}$, $h_{x_j}$ for some $i,j\in\{1,2,3\}$.
The space $\mathcal{V}_H$ corresponds to a coarsening of $\mathcal{V}_h$ obtained by condensing degrees of freedom as shown in
Fig. \ref{fig:sec2:decomposition_h_H} c). A basis for $\mathcal{V}_H$ can be obtained be selecting the functions defined on
$\Gamma$ that assume the value $1$ at some of these condensate regions (the yellow portion in Fig. \ref{fig:sec2:decomposition_h_H})
and $0$ on the rest of $\Gamma$. This basis parameterized by $\bar{H}_{x_i}$, $i=1\ldots 3$, is denoted hereafter by $\{\phi_k\}$ and
its cardinal is equal to
\begin{equation}
	\#\{\phi_j\}=
	\overbrace{\underbrace{(N_{d_1}-1)N_{d_2}N_{d_3}}_\text{\#interfaces}\underbrace{\frac{H_{x_2}}{\bar{H}_{x_2}}\frac{H_{x_3}}{\bar{H}_{x_3}}}_{\text{dof per interface}}}^\text{dof  interfaces parallel to $x_2$-$x_3$}+
	(N_{d_2}-1)N_{d_1}N_{d_3}\frac{H_{x_1}}{\bar{H}_{x_1}}\frac{H_{x_3}}{\bar{H}_{x_3}}+
	(N_{d_3}-1)N_{d_1}N_{d_2}\frac{H_{x_1}}{\bar{H}_{x_1}}\frac{H_{x_2}}{\bar{H}_{x_2}}.\label{eq:sec2:dof_PH}
\end{equation}
where $N_{d_i}=1/H_{x_i}$, $i=1,2,3$. An example of the condensation of degrees of freedom on $\Gamma$ is illustrated
in Fig. \ref{fig:sec2:decomposition_h_H} where $\bar{H}_{x_1}=2h_{x_1}$ and $\bar{H}_{x_2}=3h_{x_2}$.
There, six degrees of freedom from the fine grid, as shown in Fig. \ref{fig:sec2:decomposition_h_H} b), are condensate
into a single one, as displayed in Fig. \ref{fig:sec2:decomposition_h_H} c).

Next, fix a subdomain $\Omega^\ell$. Let us denote by $\mathcal{I}^\ell$ the indices of the functions $\phi_k$ with support contained
in $\partial\Omega^\ell$, i.e., $\mathcal{I}^\ell=\{k\in\mathbb{N}:\text{supp}(\phi_k)\subset\partial \Omega^\ell\cap\Gamma\neq \emptyset\}$.
The MBFs are computed by means of a discretization analogous to
that used for the fine grid solution, in this case subject to Robin boundary conditions.
Denoting by ${\mathbf{\hat{p}}_{\ell,k;P}}$ the column vector of cell pressure values corresponding to $\mathbf{\hat{p}}^{\ell}_h(\phi_k,0)$
and by ${\mathbf{\hat{p}}_{\ell,k;U}}$ the corresponding one for $\mathbf{\hat{p}}^\ell_h(0,\phi_k)$, the following local (fine grid)
systems must be solved
\begin{equation}
\mathcal{A}^{H}_{\ell}\,{\mathbf{\hat{p}}_{\ell,k;P}}=\mathbf{f}^{H}_{\ell,k;P},\qquad \mathcal{A}^{H}_{\ell}\,{\mathbf{\hat{p}}_{\ell,k;U}}=\mathbf{f}^{H}_{\ell,k;U}\qquad\text{for }k\in\mathcal{I}^\ell
\label{eq:sec2:multiscale_functions}.
\end{equation}
Abbreviating by $\gamma = (g_p,g_u,f)$ the forcing data and denoting by
${\mathbf{\bar{p}}_{\ell;\gamma}}$ the discrete solution of the particular problem \eqref{eq:sec2:localomegapar}
and by $\mathbf{f}^{H}_{\ell;\gamma}$ its right hand side, this problem is written as
\begin{equation}
\mathcal{A}^{H}_{\ell}\,{\mathbf{\bar{p}}_{\ell;\gamma}}=\mathbf{f}^{H}_{\ell;\gamma}~.
\label{eq:sec2:particular_problem_discrete}
\end{equation}
Now, for a subdomain totally immersed in $\Omega$ (i.e., $\Omega^\ell\cap \partial\Omega =\emptyset$)
the number of linear problems written in \eqref{eq:sec2:multiscale_functions} corresponds to 
 $$2\cdot \left|\mathcal{I}^\ell\right|=2\cdot
\left(\overbrace{2\,\frac{H_{x_2}}{\bar{H}_{x_2}}\frac{H_{x_3}}{\bar{H}_{x_3}}}^{\text{dof  interfaces $||$ to $x_2$-$x_3$}}
+\overbrace{2\,\frac{H_{x_1}}{\bar{H}_{x_1}}\frac{H_{x_3}}{\bar{H}_{x_3}}}^{\text{dof  interfaces $||$ to $x_1$-$x_3$}}
+\overbrace{2\,\frac{H_{x_1}}{\bar{H}_{x_1}}\frac{H_{x_2}}{\bar{H}_{x_2}}}^{\text{dof  interfaces $||$ to $x_1$-$x_2$}}\right),$$
in cases where each $H_{x_i}$ is equal to a constant $H$ and each $\bar{H}_{x_i}$ is equal to a constant $\bar{H}$, one has that $2\cdot \left|\mathcal{I}^\ell\right|=12\,H/\bar{H}$.
Therefore, for most of the situations considered in this work, where $\bar{H}=H$, we need to solve 12 linear problems
to construct the multiscale basis functions ${\mathbf{\hat{p}}_{\ell,k;P}}$ and ${\mathbf{\hat{p}}_{\ell,k;U}}$
for each interior subdomain $\Omega^\ell$. Typical sizes of these local problems involve a few thousand unknowns
that suggest that the use of direct solvers might be a computationally effective choice.

\subsection{The fine grid transport solver}

The classical Implicit Pressure Explicit Saturation (IMPES) method was used to perform two-phase flow simulations. In such methodology, one computes first the velocity field $\mathbf{{u}}$ by calling a Darcy solver (here, Fine Grid or MRCM) and then updates the saturation field $\mathbf{S}$. These two steps are
repeated sequentially. For the saturation step, equation \eqref{eq:sec2:Buckley-Leverett} is discretized by means of a Finite Volume Method, where the computational cells have sizes $h_{x_1}$, $h_  {x_2}$ and $h_{x_3}$ on each axis. Thus, knowing the discrete field $\mathbf{S}^n$ at time $t^n$, the update $\mathbf{S}^{n+1}=\mathbf{S}(t^n+\Delta t_s)$ is computed according to:
\begin{align}
\mathbf{S}^{n+1}_{i,j,k}=\mathbf{S}^{n}_{i,j,k}
\,&+
\frac{\Delta t_s}{h_{x_1}}\left(\varphi_{i+\frac{1}{2},j,k}\,{u}_{i+\frac{1}{2},j,k}^n+\varphi_{i-\frac{1}{2},j,k}\,{u}_{i-\frac{1}{2},j,k}^n\right)\label{eq:transporte_discreto}\\
&+\frac{\Delta t_s}{h_{x_2}}\left(\varphi_{i,j+\frac{1}{2},k}\,{u}_{i,j+\frac{1}{2},k}^n+\varphi_{i,j-\frac{1}{2},k}\,{u}^n_{i,j-\frac{1}{2},k}\right)\nonumber\\
&+\frac{\Delta t_s}{h_{x_3}}\left(\varphi_{i,j,k+\frac{1}{2}}\,{u}_{i,j,k+\frac{1}{2}}^n+\varphi_{i,j,k-\frac{1}{2}}\,{u}^n_{i,j,k-\frac{1}{2}}\right)-\nonumber\\
&-\Delta t_s\,\varphi_{i,j,k}\,\mathbf{D}\,u^n_{i,j,k}\nonumber~ ,
\end{align}
where ${u}_{i\pm\frac{1}{2},j,k}^n$, ${u}_{i,j\pm\frac{1}{2},k}^n$, ${u}_{i,j,k\pm\frac{1}{2}}^n$ denote the normal fluxes $\mathbf{u}^n\cdot\hat{\bf n}$ at the cells interfaces, and $\varphi_{i\pm\frac{1}{2},j,k}$, $\varphi_{i,j\pm\frac{1}{2},k}$, $\varphi_{i,j,k\pm\frac{1}{2}}$ are upwind approximations of $\varphi(S)$:

\begin{equation}
\varphi_{i\pm\frac{1}{2},j,k}=
\begin{cases}
\varphi(\mathbf{S}^{n}_{i,j,k})  &\text{if}~~ {u}_{i\pm\frac{1}{2},j,k} \geq 0\\
\varphi(\mathbf{S}^{n}_{i\pm 1,j,k})&\text{otherwise}
\end{cases},\label{eq:flux_disc_up}
\end{equation}
with analogous definitions for $\varphi_{i,j\pm\frac{1}{2},k}$ and $\varphi_{i,j,k\pm\frac{1}{2}}$. The term $(\Delta t_s\,\varphi_{i,j,k}\,\mathbf{D}u_{i,j,k}^n)$, with
$$\mathbf{D}u_{i,j,k}^n=u_{i+1,j,k}^n+u_{i-1,j,k}^n+u_{i,j+1,k}^n+u_{i,j-1,k}^n+u_{i,j,k+1}^n+u_{i,j,k-1}^n,$$
is added to compensate the small numerical divergence introduced by the Darcy's solver, an issue that is
particularly important when using iterative methods to obtain the solution.
This way, we assure the monotonocity of the numerical scheme.
Let us introduce a second time step $\Delta t_p$ and a parameter $C\in\{1,2,\ldots\}$ (the skipping constant)
such that $\Delta t_p\equiv C\,\Delta t_s$ along a simulation. As shown in \cite{Chen2004}, instead of solving for the
fluxes $\mathbf{u}$ at each time step $t^n$, one can freeze the field $\mathbf{u}$ and update it only every $C$ fine
time steps (with the new saturation field $\mathbf{S}(t^n+\Delta t_p)$), thus reducing the computational cost of a given
simulation while maintaining good accuracy. More details on this methodology are given below along with an overview of our implementation.

\subsection{Implementation aspects}
\label{sec:code_overview}

The implemented code has two main components: 1) The Darcy (elliptic) solver and 2) The transport (hyperbolic) solver.
A general pseudo code is shown in Algorithm \ref{alg:overview-code}. The Darcy solvers correspond to
the implementation of the numerical methods presented in Section \ref{sec:elliptic-solvers}. For both implementations
the Portable, Extensible Toolkit for Scientific Computation (PETSc)\footnote{\url{https://www.mcs.anl.gov/petsc}}, version 3.13.0, was used as interface with HPC
libraries to solve and precondition the related linear problems.

\definecolor{gray}{rgb}{0.6, 0.6,0.6}

\begin{algorithm}[H]
	\SetAlgoNoLine
	\DontPrintSemicolon
	\SetKwFunction{FRecurs}{FnRecursive}%
	\KwIn{$\mathbf{S}^0$ (initial saturation), $C$ (skipping constant), $\Delta t_s$, $T$, $solver\_type$}
	\Begin{
         {{\it \small \color{gray} $|$ Loop over time steps}} \\ 
		$t=0$\;
		$n_s= 0$, $n_p=0$\;
		\While{$(t\leq T)$}{
                ~\\
			{{\it \small \color{gray} $|$ Elliptic solver}}\\
			\If{{\it $(\mbox{mod}(n_s,C)==0)$}}{
				${\mathbf{u}_{\ell}^{n_p}}=\text{\bf DarcySolver}\left(solver\_type,\,\mathbf{S}^{n_s},{\mathbf{K}}\right)$\;
				$n_p=n_p+1$\;		
			}
                        ~\\
			{{\it \small \color{gray} $|$ Hyperbolic (fine grid) solver}}\\
			$\mathbf{S}^{{n_s}+1}=\text{\bf TransportSolver}\left(\mathbf{S}^{n_s},\,{\mathbf{u}^{n_p}_{\ell}}\right)$ \;
                        ~\\
			$t= t+\Delta t_s$\;
			$n_s= n_s+1$\;
			$\mathbf{S}^{n_s}= \mathbf{S}^{n_s+1}$\;
		}
		\caption{General overview of the in-house developed code used in this work. The variable $solver\_type$ can be set to the Fine grid or
                the Multiscale solver.}\label{alg:overview-code}
	}
\end{algorithm}

\subsubsection{Implementation details for the global fine grid solver}

\begin{figure}[ht]
	\centering
	\def\svgwidth{\linewidth}
	\includegraphics[width=0.5\linewidth]{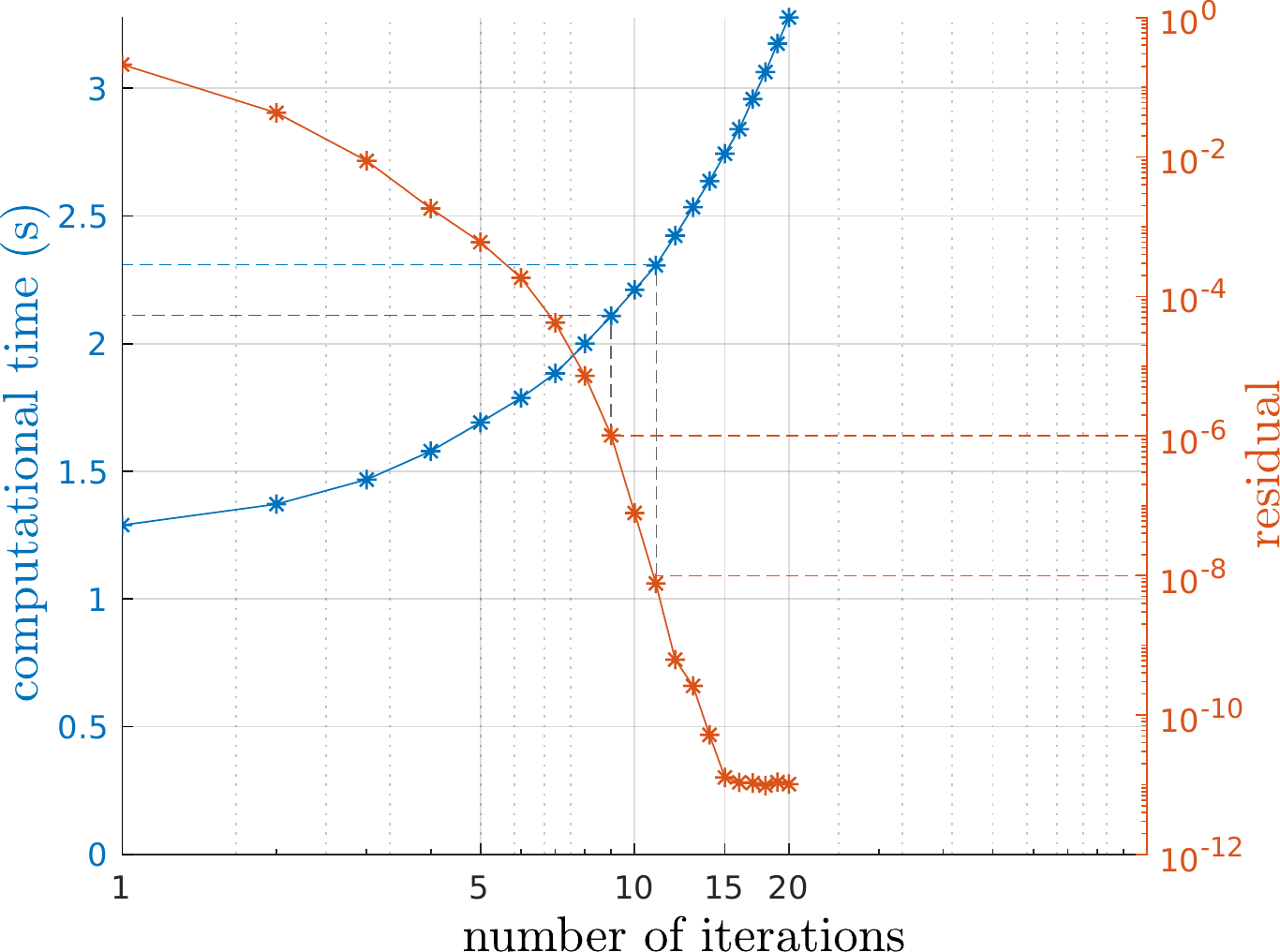}
	\caption{GMRES computational time and the solution residual as function of the number of iterations when preconditioning
        with BoomerAMG a $1$M problem (SPE10 benchmark) with 8 computational cores.}
	\label{fig:sec2:comptime-iters-gmres}
\end{figure}

A key point in this work is to use an efficient fine grid global solver in order to make a fair comparison with the proposed multiscale solver.
To that end, the linear problem \eqref{eq:sec2-darcy-discrete-fine-grid} can be solved in parallel by using state-of-the-art linear
algebra tools. In this work we choose the well known PETSc toolkit that allows to handle several combinations of linear solvers and preconditioners. Considering the sizes of the target problems we intend to consider (around 1 billion cells),
three iterative linear solvers were tested, the Generalized Minimal Residual (GMRES), the Conjugate Gradient (CG) and the BiCGStab (BCGSL) method.
Two types of preconditioners were also considered, algebraic multigrid and incomplete factorization, being the former the only one that assured convergence
and the best computational times. Among the algebraic multigrid methods available in PETSc, the best performance resulted from using either the ``Multi Level Preconditioning Package''
(ML) by Trilinos or the ``BoomerAMG'' implementation included in the library HYPRE\footnote{\url{https://computing.llnl.gov/projects/hypre-scalable-linear-solvers-multigrid-methods}}
by the Lawrence Livermore National Laboratory\footnote{\url{https://trilinos.github.io}}, both giving comparable results. However, up to its current version,
ML was limited to 32-bits indices so it was not possible to use it in high resolution problems. By fixing the BoomerAMG as the preconditioner no relevant difference on computational
times was observed when using GMRES, CG and or BCGSL. Hereafter, the results for the Fine Grid method will be relative to the GMRES solver preconditioned with BoomerAMG.
The command line options used to set that configuration are shown in Table \ref{tab:fine_grid_petsc_options}.
The option $\verb|-ksp_rtol|$ sets the tolerance for the stopping criteria based on the relative preconditioned residual.
Setting this option to \verb|1e-8| was enough to observe convergence. An increase on this tolerance only provides small saving in computing time.
To illustrate this we solve the SPE10 benchmark problem that have 
$\sim1$M unknowns distributed over 8 computational cores.
Figure \ref{fig:sec2:comptime-iters-gmres} shows the computational time spent by the Fine Grid solver and the resulting
residual as a function of the number of GMRES iterations.
A first eye catching conclusion is that the computational time slightly increases when the solution residual moves from $10^{-6}$ to $10^{-8}$.
It is also notable that the first iteration of GMRES takes most part of the time, which is related to the cost of computing the BoomerAMG operator.
Similar conclusions apply to larger problems.



\begin{table}[h]
	\begin{center}
		\begin{tabular}{cccc}
			\toprule
		$\verb|-ksp_type|$	& $\verb|-ksp_rtol|$ & $\verb|-pc_type|$ & $\verb|-pc_hypre_type|$ \\ \midrule
		$\verb|gmres|$  & $\verb|1e-8|$ & $\verb|hypre|$ & $\verb|boomeramg|$\\
			\bottomrule
		\end{tabular}
		\caption{PETSc options used to implement the Fine Grid Solver.}\label{tab:fine_grid_petsc_options}
	\end{center}
\end{table}

To close this discussion, we remark that several of options are available through the PETSc interface to configure BoomerAMG.
Preliminary tests were made by changing: the type of the multigrid cycle (\verb|V| or \verb|W| type), the
number of cycles per GMRES iteration, the type of smoother to be applied on each grid, the type of interpolator
and the use of nodal systems coarsening \cite{hypremanual}. Each one of these options has a default setting that
can be accessed by including the argument \verb|-help| in the command line. In our experience changing these
default parameters did not affect the solver's efficiency significantly, so they are the ones
adopted in our studies reported in Section \ref{sec:hpc-results-darcy}.

\subsubsection{Implementation details of the Multiscale elliptic solver}\label{sec:msc_solver_detail}

\begin{figure}[ht]
	\centering
	\def\svgwidth{\linewidth}
	\includegraphics[width=0.9\linewidth]{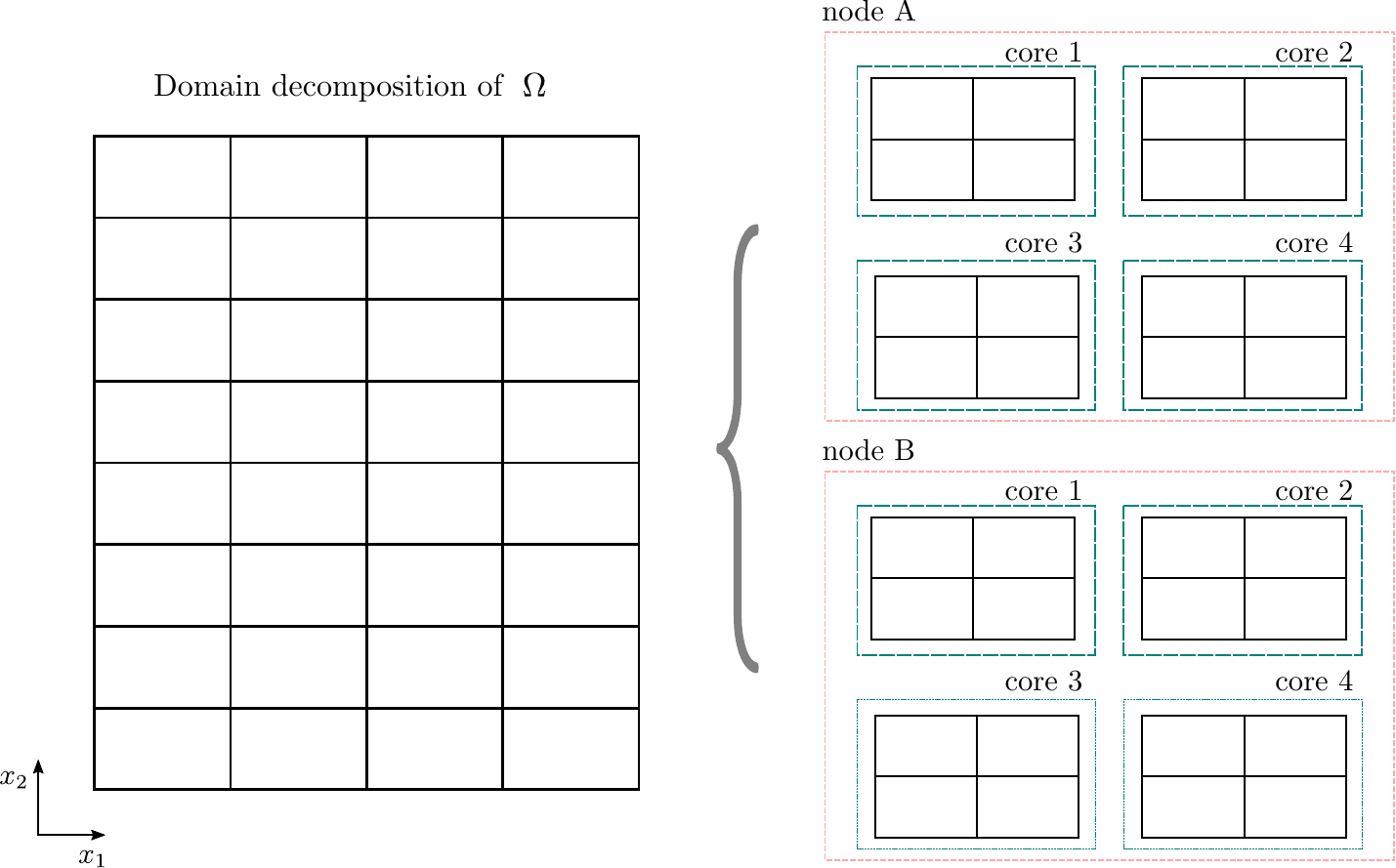}
	\caption{Example of a distribution for a decomposition of $\Omega$ into $4\times 8=32$ subdomains. The MPI processes
        are evenly spread across a cluster of 2 nodes, with 4 physical cores per node, which results in an assignment of 4 subdomains per core.}
	\label{fig:sec2:processes-distribution}
\end{figure}

Given a domain decomposition, the number of subdomains can be equal to or greater than the number of cores allocated for the MPI call (\verb|mpirun -np|).
An example of a distribution of subdomains is shown in Fig. \ref{fig:sec2:processes-distribution}. A constant number
of subdomains is distributed to a each MPI process, and hyper-threading is avoided forcing each one of these processes to run on a single physical core. A pseudo-code of the implementation is presented
in Algorithm \ref{alg:MRCM}. Note that it is composed of three \emph{loops},
all ranging over the subdomain decomposition, namely,
1) Local problems defining the MBFs are solved,
2) Interface problem is assembled and solved;
3) Multiscale solution is recovered.
Between the first two loops the interface problem is solved and so MPI communications between processors are necessary.
Note also that in order to obtain the multiscale solution, these two are the sole lines where communications are executed.
In the postprocessing step MPI calls may also be involved. Computational time profiles for each stage of the implementation
are reported in Section \ref{sec:hpc-results-darcy} for both MRCM and fine grid solvers.

\begin{algorithm}[H]
	\SetKwFunction{FRecurs}{FnRecursive}%
	\KwIn{$\tilde{\mathbf{K}}=\mathbf{K}$, processor index $r$}
        \LinesNumbered
	\SetAlgoNoLine
	\DontPrintSemicolon
	\Begin{
		{{\it \small \color{gray} $|$ Compute MBFs}} \\
		\label{line:MRCM-loop1INICIO}
		\For{$\ell \in \mathcal{S}^{r}$}{
                           
			{\bf Assembly} $\mathcal{A}^{H}_{\ell}$, $\mathbf{f}^{H}_{\ell;\gamma}$, and $\mathbf{f}^{H}_{\ell,k;P},\mathbf{f}^{H}_{\ell,k;U}$ {\bf for each} $k\in\mathcal{I}^\ell$ \;
			$\left[{\mathbf{\hat{p}}_{\ell,k;P}},{\mathbf{\hat{p}}_{\ell,k;U}},\mathbf{\bar{p}}_{\ell;\gamma}\right]=
			\text{\bf LocalSolver}\left(\mathcal{A}^{H}_{\ell},\,
			\mathbf{f}^{H}_{\ell,k;P},\,
			\mathbf{f}^{H}_{\ell,k;U},\,
			\mathbf{f}^{H}_{\ell;\gamma}\right)$ \;
                        
		}\label{line:MRCM-loop1FIM}

                ~ \\


                {{\it \small \color{gray} $|$ Interface prolem}} \\
		
		$\mathcal{A}^{\text{int}}$, $\mathbf{b}^{\text{int}} = {\bf Assembly}(
                {\mathbf{\hat{p}}_{\ell,k;P}},\,
			{\mathbf{\hat{p}}_{\ell,k;U}},\,
			\mathbf{\bar{p}}_{\ell;\gamma}) $ \;\label{line:MRCM-int1}

		$X=\text{\bf InterfaceSolver}\left(\mathcal{A}^{\text{int}},\mathbf{b}^{\text{int}}\right)$ \;\label{line:MRCM-int2}
                ~ \\

                {{\it \small \color{gray} $|$ Recover Multiscale solution}} \\
		\For{$\ell \in \mathcal{S}^{r}$ } { 
                \label{line:MRCM-loop2INICIO}    
			$\left[{{p}_{\ell}},{\mathbf{u}_{\ell}}\right]=
			\text{\bf ReconstructMscSolution}\left(\mathcal{A}^{H}_{\ell},\,
			{\mathbf{\hat{p}}_{\ell,k;P}},\,
			{\mathbf{\hat{p}}_{\ell,k;U}},\,
			\mathbf{\bar{p}}_{\ell;\gamma},X\right)$ \;
		}\label{line:MRCM-loop2FIM}
		~\\
                
		\caption{Pseudocode for the DD Multiscale mixed method.}\label{alg:MRCM}
	}
\end{algorithm}

\medskip

\subsubsection*{$|$ MBFs: {\bf LocalSolver()}}     

Let us now turn the attention to the computation of the MBFs.
Consider a subdomain $\Omega_\ell$. The computations of the MBFs and the particular problem require
solving several linear systems with the same underlying matrix $\mathcal{A}^{H}_{\ell}$.
Linear solvers based on matrix factorization are adequate for handling multiple right-hand-sides
for the typical sizes of $\mathcal{A}^{H}_{\ell}$, as anticipated.
To support this choice, let us compare the computational times obtained when solving a Darcy
problem by performing a LDL$^\intercal$ factorization
versus the computational cost corresponding to a GMRES solver with the Algebraic-Multigrid (AMG)
preconditioning \cite{Ruge1987}.
For these numerical experiments a discrete domain consisting of $45$M cells is
decomposed into 128, 256, 512, 1024, and 2048 subdomains,
such that the number of cells of each $\Omega_\ell$ ranges from $352$K  to $22$K.
Also, setting the ratio $H/\bar{H}$ to 1 and 2 ($H$-refinement in one direction),
the number of right-hand-sides ($\textbf{nrhs}$) ends up being 13 and 25, respectively.
The computational times measured when distributing the subdomains along 32 cores are reported in
Table \ref{tab:comp_local_mumps_hypre}. One observes that, even when the number of subdomains
per core is increasing, a reduction of the number of cells on each $\Omega_\ell$ results in smaller
computational times irrespective of the solver.

\begin{table}[h!]
	\begin{center}
		\begin{tabular}{r|ccccc}
			\toprule
			\#subdomains per core [{total}] & 4 [{\small 128}]  & 8 [{\small 256}] & 16 [{\small 512}] & 32 [{\small 1024}] & 64 [{\small 2048}] \\
			\midrule
			$N_{\Omega^\ell}$ & 352K & 176K & 88K & 44K & 22K \\ 
			\midrule
			LDL$^\intercal$ factorization ({\small $\mathbf{nrhs}=13$}) & 73.3 & 36.6 &20.8  & 19.1 &\textbf{14.2} \\
			GMRES + AMG ({\small $\mathbf{nrhs}=13$})& 55.3 & 35.4 & 26.2 & 24.8 & 22.1 \\ 
			\midrule
			LDL$^\intercal$ factorization ({\small $\mathbf{nrhs}=25$}) &
			74.8 & 37.9 & 21.6 & 20.0 & 14.7 \\
			GMRES + AMG ({\small $\mathbf{nrhs}=25$}) & 
			77.9 & 53.5 & 41.1 & 39.4 & 35.8 \\
			\bottomrule
		\end{tabular}
		\caption{Computing times (in seconds) for solving linear problems ($\mathcal{A}^{H}_{\ell}x=B$)
                with multiple right hand sides with two different solvers: a direct solver based on LDL$^\intercal$
                factorization, and an iterative GMRES preconditioned with an Algebraic Multigrid (AMG) method. 
                The number of righ-hand-sides is denoted by $\mathbf{nrhs}$.
                }\label{tab:comp_local_mumps_hypre}
	\end{center}
\end{table}

{\color{black}
At this point, we may ask whether it is convenient or not
to continue decreasing $N_{\Omega_\ell}$. Assuming the computational cost for solving the MBFs
for a single subdomain is of order $N_{\Omega_\ell}^p$, a simple computation shows that for solving all
the subdomains there is an advantage only if $p>1$. We experimentally find $p$
in Figure \ref{fig:sec2:times-mumps} that reports the computational times obtained by means
of the MUMPS library \cite{mumps1,mumps2} when solving a 3-dimensional Laplacian along with thirteen right-hand-sides.
For the MBFs computations, the best results are obtained by using the \verb+amd+ ordering scheme
available in MUMPS. One can notice that $p$ varies smoothly from $p\approx1$
to $p\approx 2$. Hence, if all the problems are solve sequentially
by one core, there will exist a certain lower bound $N_{\Omega_\ell}$ (typically
around $\sim 20$K for the problems of interest). Below this number the computational time associated to the local problems (MBFs)
will approach a plateau. On the other hand, moving to very large subdomains would
eventually lead to excessive computational times.
This observation will be critical later on in selecting domain decompositions that balance
the cost of the MBFs computations and the cost of the interface coupling
problem. This will be explained next.
  
  \begin{figure}[ht!]
  	\centering
  	\def\svgwidth{\linewidth}
  	\includegraphics[width=0.5\linewidth]{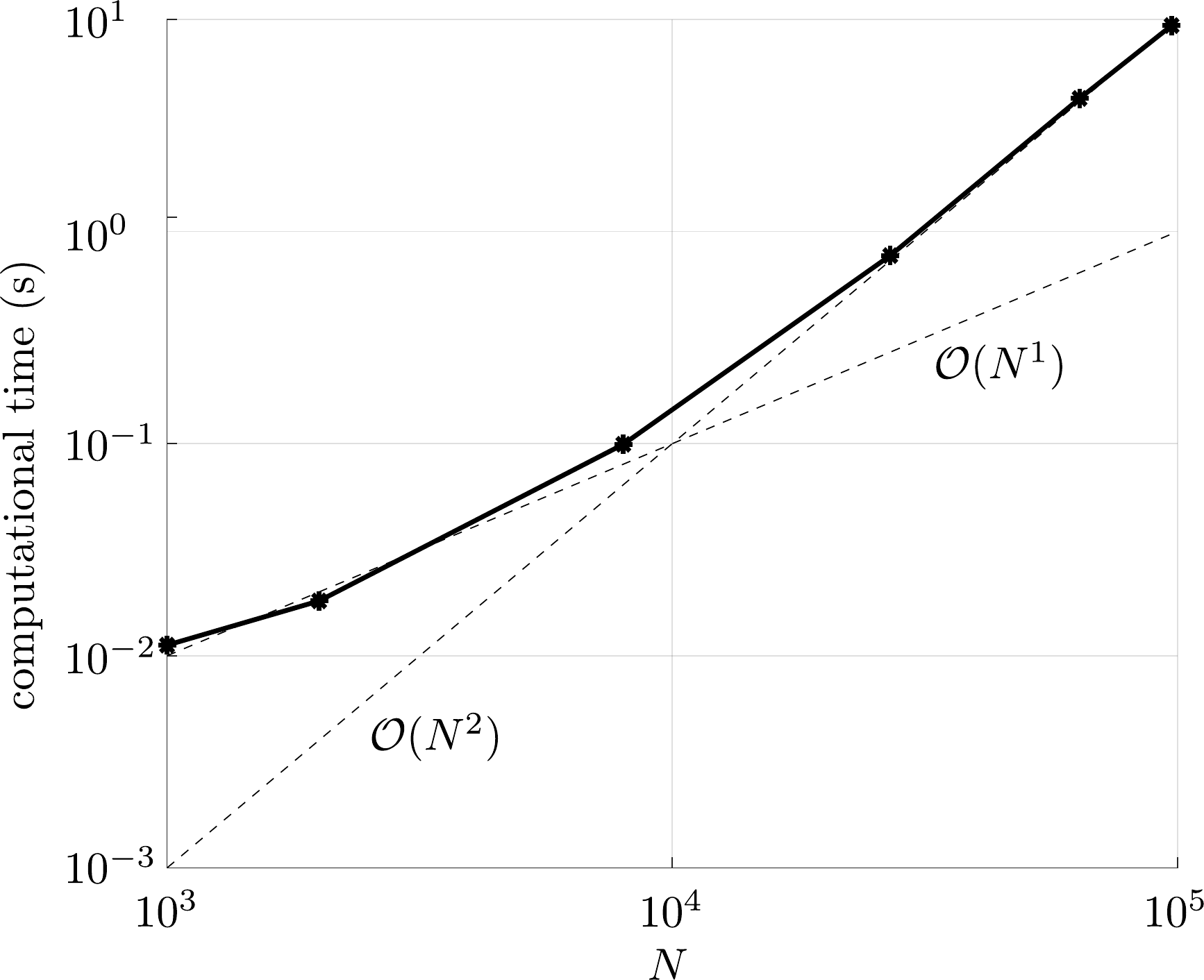}
  	\caption{Computational times versus the system order ($N$) for solving thirteen
        linear problems along LDL$^\intercal$ factorization via MUMPS.}
  	\label{fig:sec2:times-mumps}
  \end{figure}
 
\subsubsection*{$|$ Coupling problem: {\bf InterfaceSolver()}}

Based on results of the previous subsection, for billion cell reservoirs typical sizes of the
interface system \eqref{eq:interface_problem} to be solved ranges between several thousand up
to less than a million unknowns.
Hence, the computational burden of the interface resolution could be high as compared
to the MBFs computations and the choice of the linear solver is critical
for the computational efficiency.
The best option we have found is the direct method based on LU factorization implemented in MUMPS
using the \verb+pord+ ordering scheme.
To illustrate this, the computational times of LU factorization for two different configurations of the MRCM for
the same Darcy problem, the first one resulting in $431$K unknowns and
the second in $722$K unknowns, are shown in Table \ref{tab:times_interface_preliminar}.
These configurations are later on used in Section \ref{sec:weak_scaling}. On the other hand, the GMRES preconditioned with the Additive Schwarz, incomplete LU factorization and Algebraic
Multigrid methods were also tested, but none of these have provide smaller  times than LU factorization.
By means of a user-defined MPI communicator available in PETSc, one has the flexibility to choose over which MPI processes the
resolution of the interface problem is distributed.
In Table \ref{tab:times_interface_preliminar} the computational times corresponding to using 1 physical core (the root rank),
half the computational capacity of a node, a whole node, and two computational nodes are reported. One observes that increasing
the number of nodes beyond a single one does not provide significant speed-up for the $431$K case whereas for the $722$K case
some non-negligible speed-up is attained. Based on these results from now on  we adopt a direct method for the solution of the interface
system distributed on a single node.
Although not reported here for the sake of brevity, for smaller systems similar conclusions apply.

\begin{table}[h!]
	\begin{center}
		\begin{tabular}{l|cccc}
			\toprule
			$N_\Gamma$ & 1 core & 10 cores (half node) & 20 cores (full node) & 40 cores (two nodes) \\
			\midrule
			431K & 55.2 & 13.5 & 10.7 & 10.3  \\
			722K & 155.7 & 35.1 & 26.0 & 22.3 \\
			\bottomrule 
		\end{tabular}
	\caption{Computational times (in seconds) when the interface system is distribuited on different number of cores.}\label{tab:times_interface_preliminar}
	\end{center}
\end{table}

To fully assess the performance of the MUMPS library for the problems of interest in this article,
we solve the interface linear system arising from several settings to be used later on in
the HPC results section. Figure \ref{fig:mumps_interface} shows the time spent by the library to solve
problems of varying sizes $N_{\Gamma}$ from $5.5$K up to $724$K interface unknowns.
This study reveals that above $\sim 200$K unknowns the behavior of the interface
solver changes from linear to quadratic in $N_{\Gamma}$. Above such size, these
results preclude the use of direct methods for solving the interface coupling system as
is indeed verified in the experiments to be shown below.


\begin{figure}[ht!]
  	\centering
  	\def\svgwidth{\linewidth}
        \scalebox{0.55}{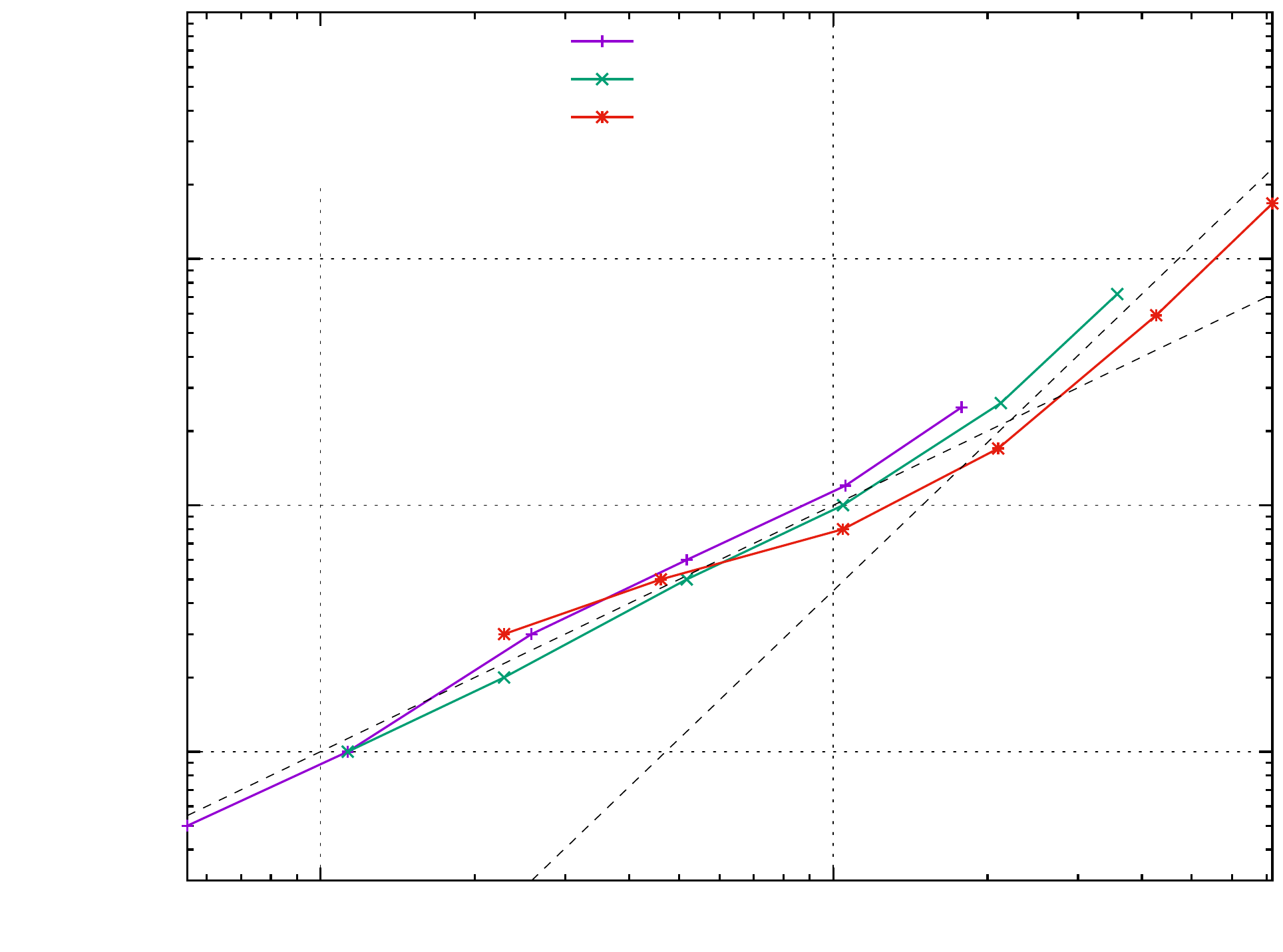}
  	\caption{Computational times versus the system order ($N_{\Gamma}$) for solving typical interface
        linear problems using $LU$ factorization via MUMPS. $P_1^{wsc}$, $P_2^{wsc}$ and
    $P_3^{wsc}$ refers to the weak scaling problems defined in Section \ref{sec:weak_scaling}.}
  	\label{fig:mumps_interface}
  \end{figure}

%% file: figs/times_mumps_interface.pdf_t
\begin{picture}(0,0)%
\includegraphics{times_mumps_interface.pdf}%
\end{picture}%
\setlength{\unitlength}{4144sp}%
\begingroup\makeatletter\ifx\SetFigFont\undefined%
\gdef\SetFigFont#1#2#3#4#5{%
  \reset@font\fontsize{#1}{#2pt}%
  \fontfamily{#3}\fontseries{#4}\fontshape{#5}%
  \selectfont}%
\fi\endgroup%
\begin{picture}(8804,6531)(903,-6753)
\put(3105,-6600){\makebox(0,0)[b]{\smash{{\SetFigFont{16}{14.4}{\familydefault}{\mddefault}{\updefault}{10$^{4}$}}}}}
\put(1500,-3281){\rotatebox{90.0}{\makebox(0,0)[b]{\smash{{\SetFigFont{16}{14.4}{\familydefault}{\mddefault}{\updefault}{computational time (s)}}}}}}
\put(5917,-6689){\makebox(0,0)[b]{\smash{{\SetFigFont{16}{14.4}{\familydefault}{\mddefault}{\updefault}{$N_{\Gamma}$}}}}}
\put(4740,-566){\makebox(0,0)[rb]{\smash{{\SetFigFont{16}{14.4}{\familydefault}{\mddefault}{\updefault}{$P_1^{wsc}$}}}}}
\put(4740,-840){\makebox(0,0)[rb]{\smash{{\SetFigFont{16}{14.4}{\familydefault}{\mddefault}{\updefault}${P_2^{wsc}}$}}}}
\put(4740,-1110){\makebox(0,0)[rb]{\smash{{\SetFigFont{16}{14.4}{\familydefault}{\mddefault}{\updefault}{$P_3^{wsc}$}}}}}
\put(6628,-6600){\makebox(0,0)[b]{\smash{{\SetFigFont{16}{14.4}{\familydefault}{\mddefault}{\updefault}{10$^{5}$}}}}}
\put(2105,-3755){\makebox(0,0)[rb]{\smash{{\SetFigFont{16}{14.4}{\familydefault}{\mddefault}{\updefault}{10$^{0}$}}}}}
\put(2105,-2062){\makebox(0,0)[rb]{\smash{{\SetFigFont{16}{14.4}{\familydefault}{\mddefault}{\updefault}{10$^{1}$}}}}}
\put(2105,-369){\makebox(0,0)[rb]{\smash{{\SetFigFont{16}{14.4}{\familydefault}{\mddefault}{\updefault}{10$^{2}$}}}}}
\put(3736,-5011){\makebox(0,0)[rb]{\smash{{\SetFigFont{16}{14.4}{\familydefault}{\mddefault}{\updefault}{$\mathcal{O}(N_{\Gamma}^1)$}}}}}
\put(6976,-4201){\makebox(0,0)[lb]{\smash{{\SetFigFont{16}{14.4}{\familydefault}{\mddefault}{\updefault}{$\mathcal{O}(N_{\Gamma}^2)$}}}}}
\put(2105,-5449){\makebox(0,0)[rb]{\smash{{\SetFigFont{16}{14.4}{\familydefault}{\mddefault}{\updefault}{10$^{-1}$}}}}}
\end{picture}%

%% file: tex/hpc_biphasic.tex
As a final numerical experiment we solve a two-phase flow problem by solving \eqref{eq:sec2:Buckley-Leverett}-\eqref{eq:sec2:bcu}
(see Algorithm \ref{alg:overview-code}). Results are obtained for both the Multiscale mixed method and the Fine Grid solver.
According to \eqref{eq:sec2:weakcont_f} the velocity solutions produced by the MRCM are conservative in a scale that corresponds
to the support of the basis functions that span the pressure interface space $\mathcal{P}_h$. As the supports are usually chosen such that
$\bar{H}_{x_j}\gg h_{x_j}$, the velocity solutions are in general discontinuous at the fine level $h_{x_j}$, except for very
large values of the algorithmic parameter $\alpha$.
A postprocessing of the velocity field is necessary prior to solving the hyperbolic transport problem. Here the \emph{Mean}
method \cite{guiraldello2020} has been chosen to that end for the sake of simplicity.
A unique flux over $\Gamma$ is defined based on the average value of the velocity on interfaces between adjacent subdomains.
These fluxes are then used as Neumann boundary conditions to compute local
problems on each $\Omega^{\ell}$. Given the multiscale solution ${\bf u}_H$, this amounts to computing
the unique velocity
\begin{equation*}
\bar{\bf U}^{\Gamma}_h \doteq \dfrac{1}{2} \left( {\bf u}^{+}_H\vert_{\Gamma} + 
{\bf u}^{-}_H\vert_{\Gamma}\right).
\end{equation*}
that in line with \eqref{eq:sec2:weakcont_f} transfers the same mass across the interfaces
as the multiscale solution. For each subdomain $\Omega_i$, we solve
\begin{equation}
\left \{
\begin{array}{rclll}
\widetilde{\bf u}^{i}_h &=& -K \, \nabla \widetilde{p}^{\,i}_h &&\  \mbox{in} \ \Omega_{i}, \\
\nabla\cdot \widetilde{\bf u}_h^{i} &=& f &&\  \mbox{in}  \ \Omega_{i}, \\
\widetilde{\bf u}^{i}_h\cdot\check{\bf n}^i &=& {\bf u}_H|_{\partial\Omega_i}\cdot\check{\bf n}^i &&\  \mbox{on}  \ \partial\Omega_{i}\cap\partial\Omega, \\
\tilde{\bf u}^{i}_h\cdot\check{\bf n}^i &=& \bar{\bf U}_h^{\Gamma}|_{\partial\Omega_i} \cdot \check{\bf n} &&\ \mbox{on} \ \partial\Omega_{i}\cap\Gamma , 
\end{array}
\right.,
\label{eq:locprob_submean}
\end{equation} 
that is conservative on the fine scale.
The procedure involves communications to compute the unique flux. For it, a single \texttt{MPI\_Allreduce} call is sufficient
having a negligible cost in the overall computational time. However, the solution of the local problems \eqref{eq:locprob_submean} involves a
new linear subdomain solve. Certainly, this is not the best method as reported in \cite{guiraldello2020} where alternative
computationally more efficient and accurate postprocessing methods can be found, although their implementation is a bit more involved
and is left for future work.

The multiscale performance is assessed by comparing the production curves. We choose the computational setting adopted for the
$P^{wsc}_2$ distributed over $288$ cores. The simulation is performed until $T_{\text{PVI}} = 0.2$, well beyond the breakthrough time for all production wells. The frequency of Darcy solves
comes from a skipping constant $C=600$. The time step $\Delta t_s$ in \eqref{eq:transporte_discreto} satisfies
the CFL condition \cite{courant1928partiellen,coats2000note}, which translates to
\begin{equation}
\Delta t_s \leq \frac{\min\left\{h_{x_1},h_{x_2},h_{x_3}\right\}}{\displaystyle\max_{\Omega}\left|\varphi'(\mathbf{S}^n)\mathbf{u}^n\right|}~.\label{eq:CFL}
\end{equation}
The oil is being extracted from the four wells located at the corners of $\Omega$.
The production ($\mathcal{P}_{\text{oil}})$ and water-cut ($W_i$) curves correspond to the fraction of oil and water
for each production well, $i=1\ldots 4$ as a function of time. The $\mathcal{P}_{\text{oil}}(t)$ is computed
according to
\begin{equation*}\label{eqn:prodcfun}
	\mathcal{P}_{\text{oil}}(t) =  \frac{\sum_{i=1}^{4}\int_{\partial\Omega_{\text{w}}^{i}}\left(1-\varphi(S_{w}(\mathbf{x},t))\right)\mathbf{u}(\mathbf{x},t)\cdot d\boldsymbol\Gamma}{\sum_{i=1}^{4}\int_{\partial\Omega_{\text{w}}^{i}}\mathbf{u}(\mathbf{x},t)\cdot d\boldsymbol\Gamma}~,
\end{equation*}
whereas for each production well ($i=1\ldots 4$) one has
\[
W_{i}(t) = \frac{\int_{\partial\Omega_{\text{w}}^{i}}\varphi(S_{w}(\mathbf{x},t))\,\mathbf{u}(\mathbf{x},t)\cdot d\boldsymbol\Gamma}{\int_{\partial\Omega_{\text{w}}^{i}}\mathbf{u}(\mathbf{x},t)\cdot d\boldsymbol\Gamma}~.
\]
The time variable used to present results is in $T_{\text{PVI}}(t)$ units defined as:
\[
T_{\text{PVI}}(t) = \frac{1}{V_p}\int_{0}^{t}\int_{\Omega_{\text{w}}^{0}} f(\mathbf{x},\tau) \,d\mathbf{x}\,d\tau,
\]
being $V_{p}$ the reservoir's total pore-volume and $f$ the source term \eqref{eq:sec2:continuity}. 

\begin{figure}[t]
\newcommand{\mw}{.49}
\begin{center}
\includegraphics[width=\mw\textwidth]{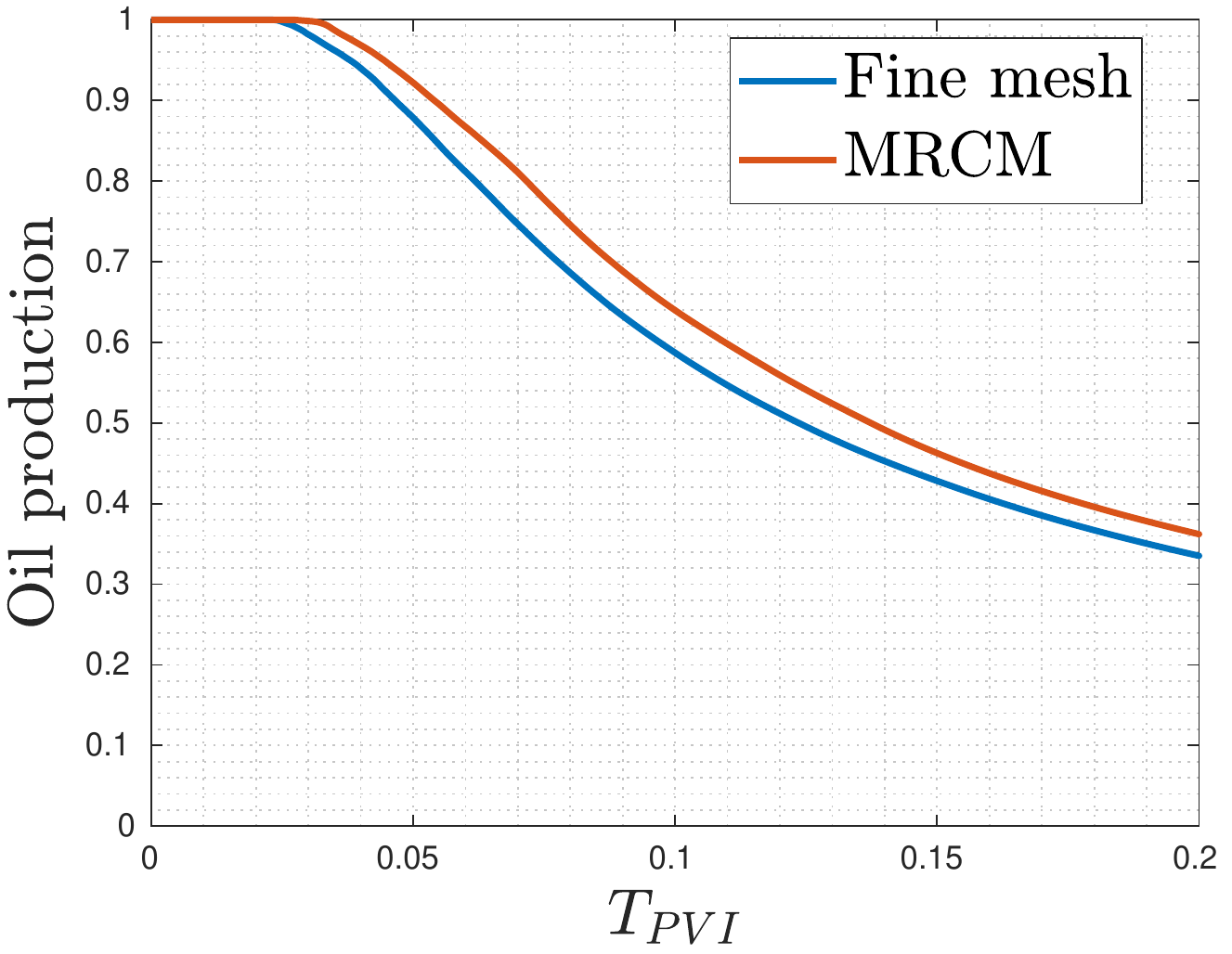}
\includegraphics[width=\mw\textwidth]{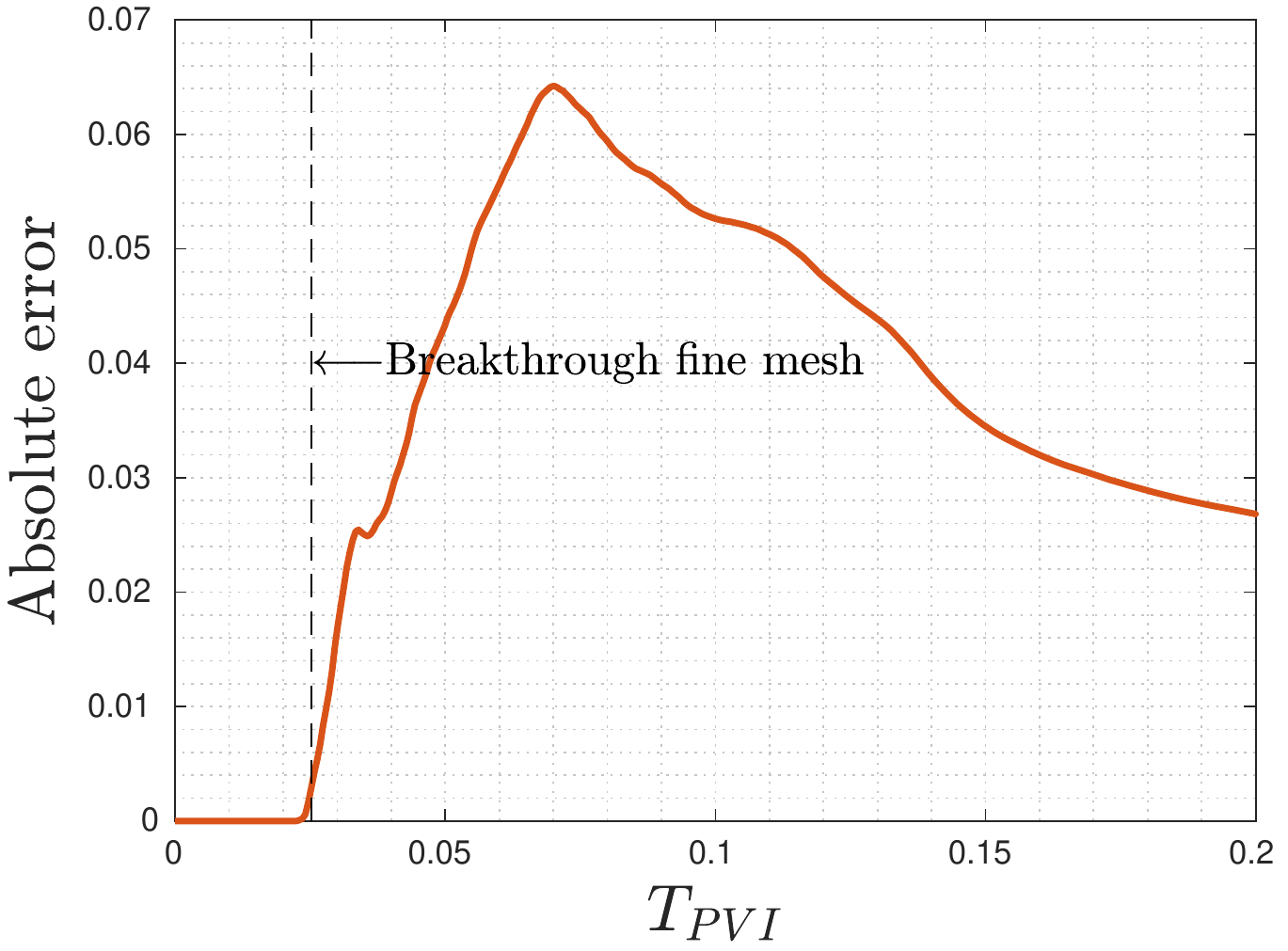}
\end{center}
\caption{Oil production curve resulting when using fine mesh and MRCM for elliptic problems (left). Absolute error between these curves (right).}\label{fig:oilperr}
\end{figure}

\begin{figure}[!t]
\newcommand{\mw}{.45}
\begin{center}
\includegraphics[width=\mw\textwidth]{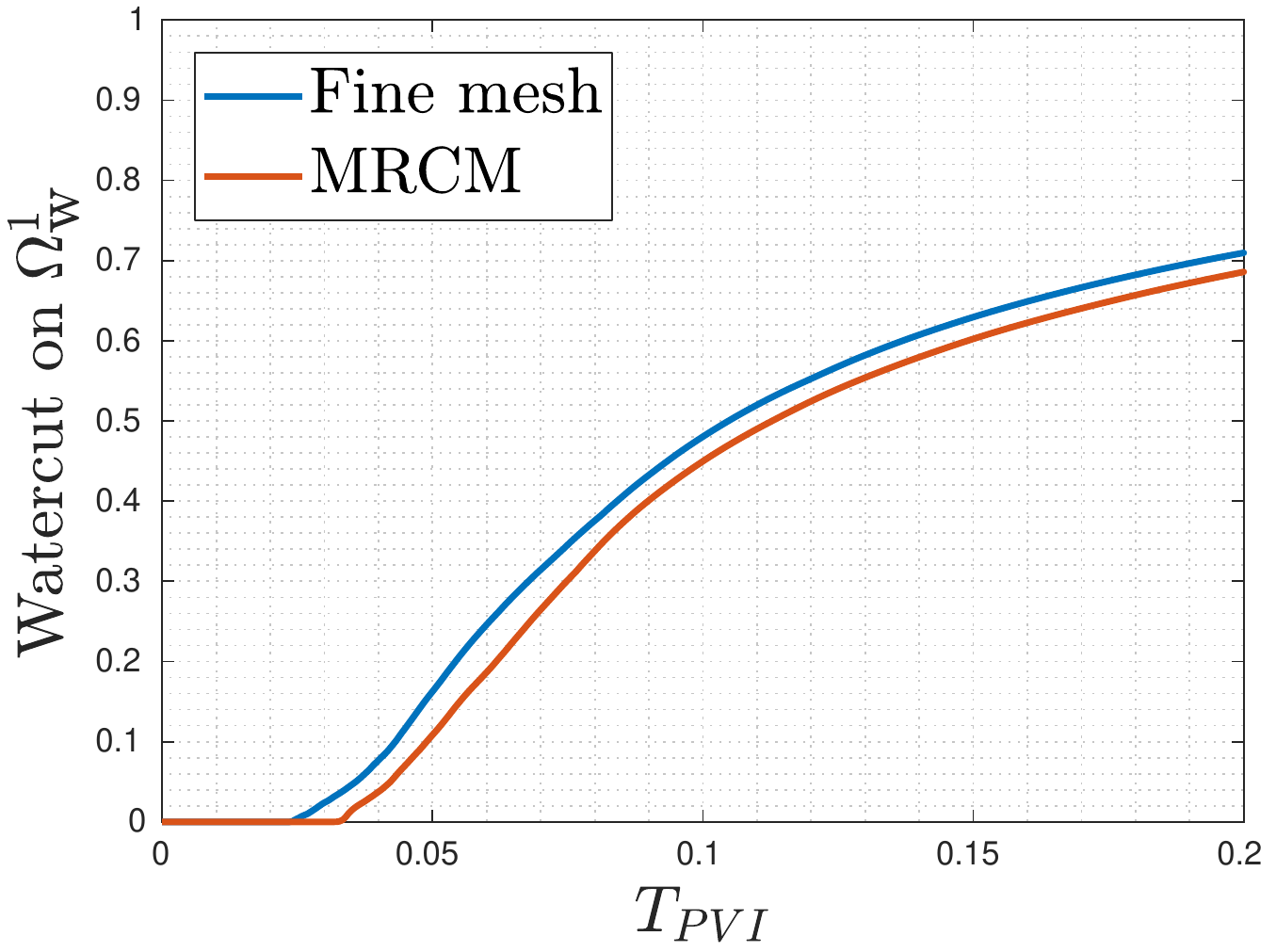}
\includegraphics[width=\mw\textwidth]{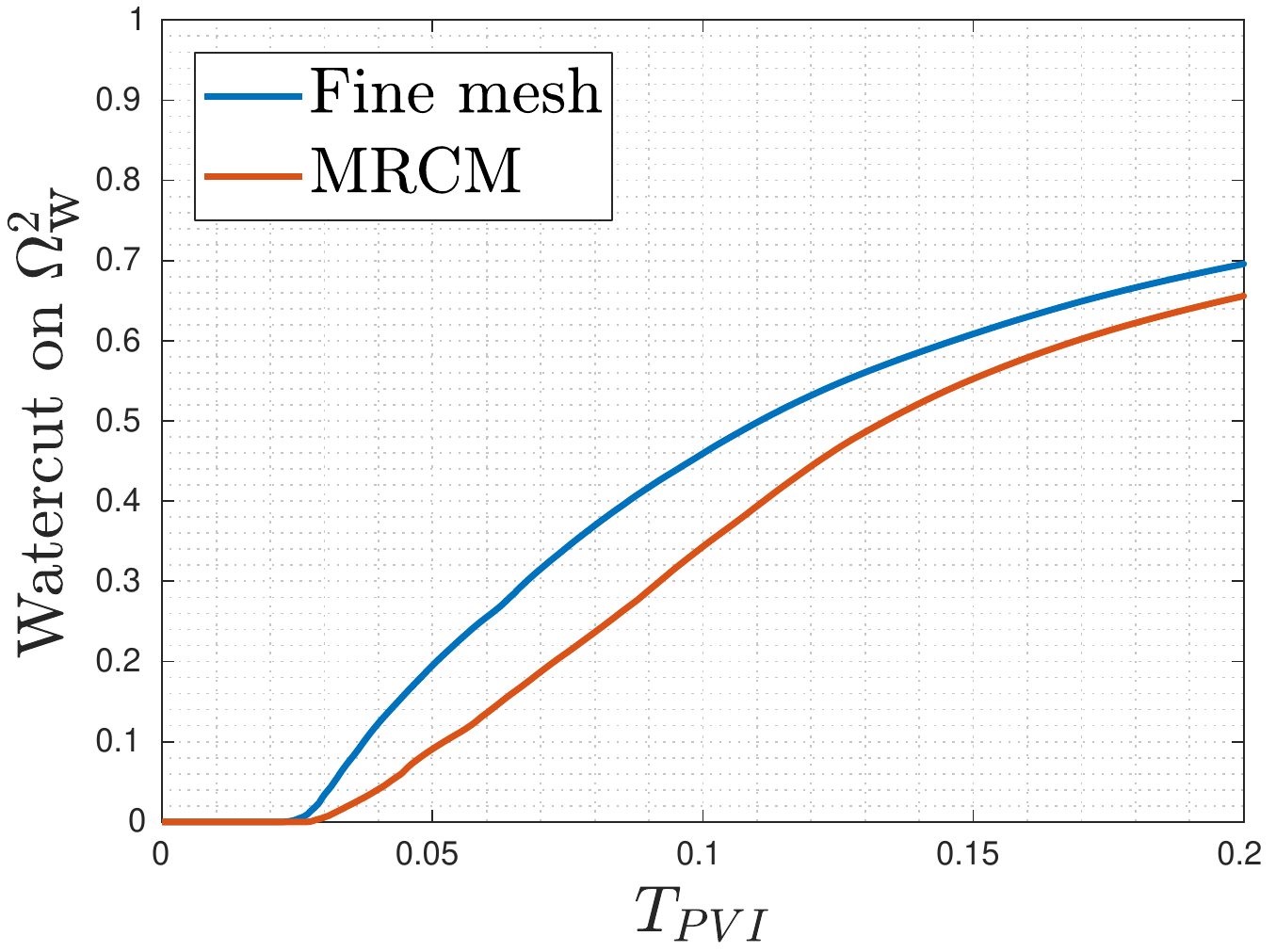}
\includegraphics[width=\mw\textwidth]{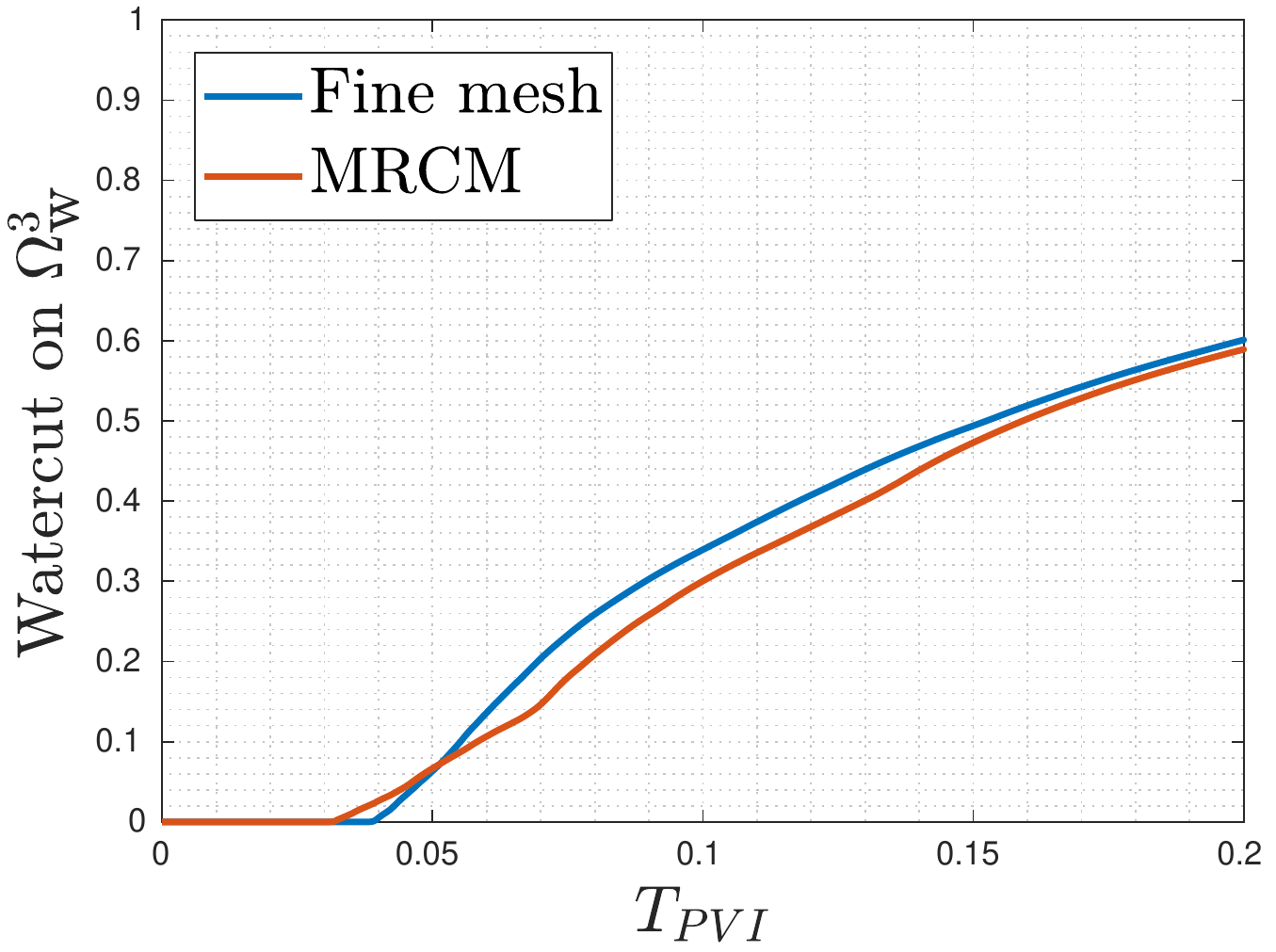}
\includegraphics[width=\mw\textwidth]{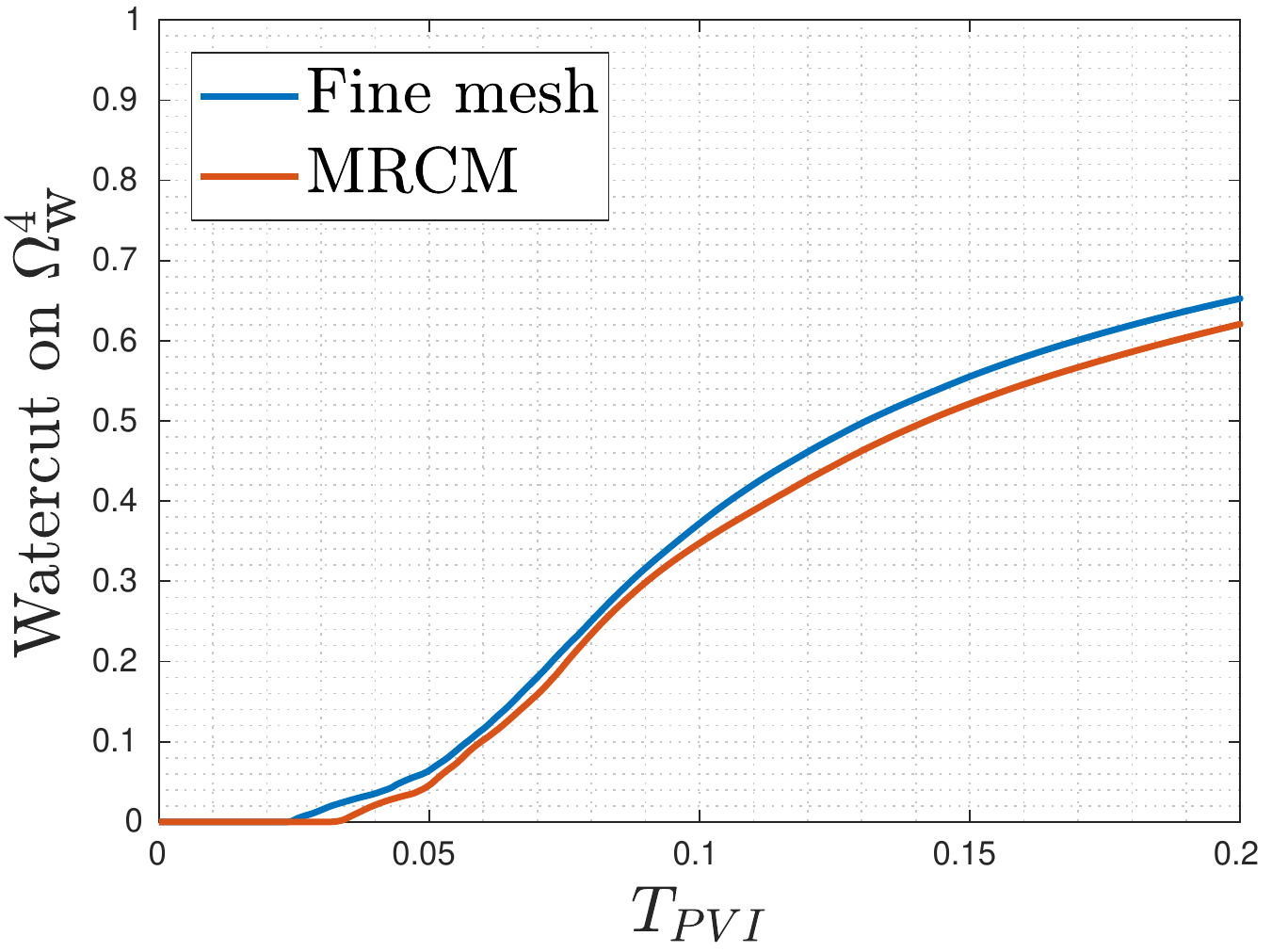}
\end{center}
\caption{Watercut on wells $\Omega_{\text{w}^{i}}$, $i = 1,2,3,4$.}\label{fig:wcwells}
\end{figure}

The computational times for the complete simulations are $473$ hours for the Fine grid solver
and $190$ hours for the multiscale method. Since the number of Darcy calls differs from one case
to the other due to adaptivity of $\Delta{t_s}$, comparison of total times is less meaningful.
The computational time of one single call being a more 
representative figure. For the considered setting the cost of a single call for the elliptic solver
when using the MRMC takes $29$ seconds, whereas it takes $83$ seconds for the fine grid solver.
Recall that prior to solving the transport equation after each Darcy solve one needs to execute
the velocity post-processing which as mentioned can be significantly reduced in the future
by implementing more efficient techniques.
To conclude, a comparison of the oil production curves resulting from the two-phase flow solver
using the Fine Grid and the MRCM is shown in Figure \ref{fig:oilperr}. There is a good agreement
of the multiscale solution to the fine grid one.
The maximum difference in the oil fraction produced is around $6\%$, which takes place at $T_{\text{PVI}} = 0.07$,
after which the difference decreases monotonically. By looking at the watercuts curves on each production well,
plotted in Figure \ref{fig:wcwells}, most part of the error is concentrated at wells $2$ and $1$.

